%% file: paper.tex
\documentclass[12pt]{article}
\usepackage{amsfonts,amsmath}
\usepackage{tikz}
\usepackage{pgfplots}
\usepackage{pgfplotstable}
\usepackage{booktabs}
\usepackage{float}
\usepackage{siunitx}

\usepackage{a4wide,epsfig}
\usepackage[utf8]{inputenc}   % damit man Umlaute direkt eingeben kann und diese erkannt werden.

\newtheorem{theorem}{Theorem}[section]

\newtheorem{corollary}[theorem]{Corollary}
\newtheorem{lemma}[theorem]{Lemma}

\newtheorem{remark}{Remark}[section]
\newcommand{\proof} [1]
   { \noindent {\bf Proof.} #1 \hfill\rule{0.5em}{1.2ex} \par\medskip}

\renewcommand{\vec}[1]{\underline{#1}}

\usepackage{hyperref}

\numberwithin{equation}{section}

\begin{document}

\setcounter{page}{1}

\title{Optimal complexity solution of space-time finite element systems for
state-based parabolic \\ distributed optimal control problems}
\author{Richard~L\"oscher, Michael~Reichelt, Olaf~Steinbach}
\date{Institut f\"ur Angewandte Mathematik, TU Graz, \\ 
Steyrergasse 30, 8010 Graz, Austria}

\maketitle

\begin{abstract}
  In this paper we consider a distributed optimal control problem subject
  to a parabolic evolution equation as constraint. The approach presented
  here is based on the variational formulation of the parabolic evolution
  equation in anisotropic Sobolev spaces, considering the control in
  $[H_{0;,0}^{1,1/2}(Q)]^*$. The state equation then defines an isomorphism
  onto $H^{1,1/2}_{0;0,}(Q)$, where the anisotropic Sobolev norm can be
  realized by using a modified Hilbert transformation. As in the classical
  approach we consider the associated optimality system from which we
  can eliminate both the adjoint state and the control. After discretization,
  the cost or regularization term involves a Schur complement matrix
  $S_h = B_h^\top A_h^{-1} B_h$, defining a norm on the state space, with
  $B_h$ being the space-time finite element matrix for the parabolic
  evolution equation, and $A_h^{-1}$ representing the norm for the control.
  While the structure of the Schur
  complement matrix $S_h$ may complicate an efficient numerical solution
  of the discrete optimality system, instead of $S_h$ we use a spectrally
  equivalent matrix $D_h$ which also defines a norm on the state space, but
  allows a more efficient realization. Moreover,
  using a space-time tensor product mesh, we can use sparse factorization
  techniques to construct a solver of almost linear complexity. 
\end{abstract}

\section{Introduction}
The analysis for optimal control problems constrained by partial differential
equations already bears a long history \cite{Lions:1971} and is by now
well-established. These problems occur in the modeling of many applications,
as, e.g., cancer treatment \cite{Deuflhard:2012,SchielaWeiser:2010}, or
shape optimization \cite{GanglKoetheCesaranoMuetze:2022}, see also
\cite{Herzog:2022} for an overview on recent advances. The overall aim
is to determine a control, which produces a related state as the unique
solution of a partial differential equation, that is as close to a desired
target as possible, while being of acceptable cost. The solution of the
partial differential equation defines the control-to-state map, and we
assume that this operator defines an isomorphism between the control and the
state space. In the standard approach this control-to-state map is used to
eliminate the state, and it remains to minimize a reduced cost functional
over the set of admissible controls. In the case of neither state nor control
constraints, the minimizer is then characterized as solution of the gradient
equation which links the control to the adjoint state as solution of the
adjoint partial differential equation. Again, when eliminating the adjoint
state, we end up with a composed operator equation to determine the control
from which we can compute the state in a post processing procedure. Instead
of eliminating the state, we may consider the equivalent optimality system
including the state, the control, and the adjoint state. Or, in contrast to
this more standard approach, we can use the partial differential equation
to eliminate the control where we now end up with a reduced cost functional
to be minimized over the set of admissible states. This in turns results in
an operator equation to determine the optimal state from which we can compute
the optimal control in a post processing procedure. On the continuous level,
all three approaches are equivalent, even when including state or control
constraints. Although both approaches may exhibit different features when
considering their numerical approximations, in general they require the
solution either of a coupled discrete block system, or of an equivalent
discrete Schur complement system. This seems to be inadequat with respect
to both the number of degrees of freedom, and to the computational complexity.
While the isomorphic control-to-state map is uniquely defined by the
underlying partial differential equation, there is some freedom in the
choice of the regularity of the control, implying the regularization norm
to be used. Note that this norm might be determined by the costs one is
interested in, when considering a particular application. In any case, since the
control-to-state map is assumed to be an isomorphism, the control
space then implies the state space, and vice versa.
When using the inverse of the control-to-state map, we can rewrite the
cost or regularization norm of the control as a norm of the related state,
the energy norm, whose realization would require the solution of the
underlying partial differential equation. However, we can use any equivalent
norm in the state space to simplify the numerical solution.

In order to balance the error between the state and the target with the
costs of the control, a cost or regularization parameter is used.
It turns out that, depending on the regularity of the target, the
error between the state and the target can be bounded in terms of this
regularization parameter, yielding convergence as this parameter
tends to zero. On the other hand, if the target is not in the state
space, which is of interest in many applications, the costs become
unbounded. The prescription of maximal costs therefore defines a lower bound
for the regularization parameter, and in the sequel, this also limits the
accuracy in the approximation of the target by the state, already on the
continuous level. When considering a discretization scheme to compute a
numerical approximation of the unknown state, this also defines a minimal
discretization parameter, e.g., the finite element mesh size, in
order to reach the same accuracy between the target and the numerical
approximation of the state. Independent of the regularity of the target
one can express this optimal mesh size in terms of the regularization
parameter. However, from a practical point of view, the regularization
parameter is chosen in terms of the mesh size as long as the costs are
acceptable. Such a strategy fits perfectly in a nested iteration scheme,
as discussed in \cite{LLSY_NLA:2023} for a distributed optimal control
problem subject to the Poisson equation.
Note that this general behavior is valid for any choice of state and
control spaces. A different choice of the control space only affects
the dependency of the regularization parameter as a function of the
mesh size. As a consequence, the resulting discrete system matrices are
symmetric, positive definite, and spectrally equivalent to diagonal
matrices, hence we can use a preconditioned conjugate gradient
method as an iterative solver.

In the case of a distributed optimal control problem subject to the
Poisson equation, using energy regularization, the reduced
optimality system turned out to be a scalar singularly perturbed
diffusion equation  \cite{LSY:2022CMAM}. The aim of this paper
is to provide a similar setting for a distributed optimal control
problem subject to the heat equation with zero Dirichlet boundary and
initial conditions in the space-time domain
$Q = \Omega \times (0,T)$. In this case, the
most standard approach is to measure the costs of the control in
$L^2$, allowing for a smooth, convex functional to be minimized, and
ensuring a unique solution of the control problem, see, e.g.,
\cite{Lions:1971,Troeltzsch:2010}. Though, in a wide variety of
applications, this approach is valid, in recent years there was a growing
interest in sparse controls \cite{Casas:2017,Stadler:2009}. Their
significant advantages are in providing information about the optimal
location of control devices, rather than only their intensity. To gain
these sparsity results, there are different approaches in the literature,
where either the $L^1$ norm of the control is added in the cost
functional \cite{CasasHerzogWachsmuth:2017,Casas:2013,
  HerzogStadlerWachsmuth:2012}, or the control is considered to be of
bounded variation \cite{CasasKunisch:2017}, or to be measure valued
\cite{CasasKunisch:2016,KunischPieperVexler:2014}, allowing, e.g.,
controls as distributions. While in the standard setting the control
is considered in $L^2(Q)$ and the state in the Bochner space
$L^2(0,T;H^1_0(\Omega))$, see, e.g., \cite{LSTY_SISC:2021}, this setting
does not define an isomorphism. In order to define an isomorphism, we can
consider $L^2(0,T;H^{-1}(\Omega))$ as control space, as it was done
in \cite{LSTY_SINUM:2021, LSY:2022}. But in all of these cases we
finally have to solve a coupled system of forward and backward
heat equations. 

Following \cite{SteinbachZank:2020} we can also consider a variational
formulation of the heat equation in anisotropic Sobolev spaces with
the state in $H^{1,1/2}_{0;0,}(Q)$, and the control in
$[H^{1,1/2}_{0;,0}(Q)]^*$. While the numerical analysis of space-time
finite element methods for this approach is similar as in the
case of energy regularization in Bochner spaces, the main difference
is in the computability of the regularization norm. Since the norm of
the control in $[H^{1,1/2}_{0;,0}(Q)]^*$ is equivalent to a norm of the
state in $H^{1,1/2}_{0;0,}(Q)$, for the latter
we can choose a norm which is simpler to compute, i.e.,
using the modefied Hilbert transformation
${\mathcal{H}}_T$ as introduced in \cite{SteinbachZank:2020} we have
\[
  \| u \|^2_{H^{1,1/2}_{0;0,}(Q)} =
  \langle \partial_t u , {\mathcal{H}}_T u \rangle_Q +
  \| \nabla_x u \|^2_{L^2(Q)} .
\]
When using this norm, which also implies a norm for the control
in $[H^{1,1/2}_{0;,0}(Q)]^*$, it remains to solve a scalar singularly
perturbed elliptic partial differential equation in the space-time domain.
Although this approach can be used for completely unstructed
discretizations in space and time, an efficient evaluation of the
modified Hilbert transformation ${\mathcal{H}}_T$, which is non-local in time,
is yet not available for this general setting. Instead
we will consider a space-time tensor product structure and apply
strategies from \cite{LangerZank:2021} to present a solver that is of
optimal complexity and fiercely easy to parallelize both in space and
time.       

The main new contribution of this paper is on an efficient implementation
of the anisotropic Sobolev norm in $H^{1,1/2}_{0;0,}(Q)$ to be used as
regularization or cost for the minimization of certain functionals in
optimal control problems. While we consider a tracking type functional
subject to the heat equation with a distributed control as a model problem,
this approach is applicable also to optimal control problems with
partial observations or controls, and to other functionals including this
anisotropic Sobolev norm. But such extensions exceed the scope of this
paper. However, we will also present some numerical results in the
case of additional state constraints which are of interest in many
applications, see also \cite{löscher2024optimalcomplexitysolutionspacetime}.

The rest of the paper is organized as follows:
In Section \ref{Section:Review} we review some existing approaches for the
space-time solution of a tracking-type distributed optimal control problem
subject to the heat equation. In particular, this involves both the
regularization in $L^2(Q)$, and in the energy space
$L^2(0,T;H^{-1}(\Omega))$. As an alternative to these approaches, in
Section \ref{Section:Aniso} we introduce a new regularization for the
control in the dual space
of the anisotropic Sobolev space $H^{1,1/2}_{0;,0}(Q)$. While we can
analyze this as in the previous cases, using that the heat equation defines an isomorphism, we can reformulate the minimization problem with
respect to the state, where the regularization for the state is now
considered in $H^{1,1/2}_{0;0,}(Q)$. Using a modified Hilbert transformation
we are able to derive a computable representation of this
anisotropic Sobolev norm. In particular, in Section \ref{Section:Solver}
we analyze the Galerkin discretization of the first order temporal
derivative in combination with the modified Hilbert transformation
which results in a symmetric and positive definite stiffness matrix
$A_{h_t}$.
In the case of a uniform discretization in time we compute all
eigenvectors and related eigenvalues of a generalized eigenvalue problem
$A_{h_t} \underline{v} = \lambda M_{h_t} \underline{v}$ with the temporal
mass matrix, using piecewise linear continuous basis functions. In the
case of a space-time tensor product discretization, this then allows 
for an efficient solution of the global problem with optimal, i.e.,
almost linear, complexity.
In Section \ref{Section:Numerics} we first present three
examples for target functions of different regularity, which confirm
all the theoretical results. In addition, and as in previous work,
we also present the turning wave example in order to demonstrate
the applicability and efficiency of the proposed approach in the
case of a non-linear state equation \cite{LSTY_SISC:2021}. Finally,
we present an example with state constraints where the computational
domain is a human gallbladder.

\section{A parabolic distributed optimal control problem}
\label{Section:Review}
As a model problem we consider the minimization of the tracking type
functional
\begin{equation}\label{functional}
  {\mathcal{J}}(u_\varrho,z_\varrho) =
  \frac{1}{2} \int_Q [u_\varrho - \overline{u}]^2 \,
  dx \, dt + \frac{1}{2} \, \varrho \, \| z_\varrho \|^2_Z
\end{equation}
subject to the Dirichlet boundary value problem for the heat equation,
\begin{equation}\label{heat equation}
  \partial_t u_\varrho - \Delta u_\varrho = z_\varrho \; \mbox{in} \; Q, \quad
  u_\varrho = 0 \; \mbox{on} \; \Sigma, \quad u_\varrho = 0 \;
  \mbox{on} \; \Sigma_0,
\end{equation}
where $Q := \Omega \times (0,T)$, $\Sigma := \partial \Omega \times (0,T)$,
and $\Sigma_0:=\Omega \times \{ 0 \}$.
Here, $\Omega \subset {\mathbb{R}}^n$, $n=2,3$, is a bounded domain with
Lipschitz boundary $\partial \Omega$, and $T>0$ is a given time horizon.
Finally, $\varrho \in {\mathbb{R}}_+$ is some regularization or cost
parameter on which the solution depends.
In \eqref{functional}, we aim to approximate a given target
$\overline{u} \in L^2(Q)$ by a function $u_\varrho$ satisfying the heat equation
\eqref{heat equation} with the control $z_\varrho$ as volume source.
Of particular interest is the case when the target $\overline{u}$ is not
in the state space, e.g., $\overline{u}$ is discontinuous, or violates
homogeneous boundary and initial conditions.

The space-time variational formulation,
e.g., \cite{SchwabStevenson:2009, Steinbach:2015, UrbanPatera:2014},
of the initial boundary value problem \eqref{heat equation} reads to find
$u_\varrho \in X_0 := \{ u \in Y : \partial_t u \in Y^*, u = 0 \;
\mbox{on} \; \Sigma_0 \}$ such that
\begin{equation}\label{VF Bochner}
  \langle \partial_t u_\varrho , v \rangle_Q +
  \langle \nabla_x u_\varrho , \nabla_x v \rangle_{L^2(Q)} =
  \langle z_\varrho , v \rangle_Q
\end{equation}
is satisfied for all $v \in Y := L^2(0,T;H^1_0(\Omega))$. Note that
$\langle z_\varrho , v \rangle_Q$ denotes the duality pairing for $v \in Y$ and
$z_\varrho \in Y^*$ as extension of the inner product in $L^2(Q)$.
The solution of the space-time variational formulation
\eqref{VF Bochner} defines the control-to-state map
$u_\varrho = {\mathcal{S}}z_\varrho$, and hence we can write the reduced
cost functional, in the particular case $Z=L^2(Q)$, as
\begin{equation}\label{reduced L2 standard}
  \widetilde{\mathcal{J}}(z_\varrho) = \frac{1}{2} \, \int_Q
  [{\mathcal{S}}z_\varrho - \overline{u}]^2 \, dx \, dt +
  \frac{1}{2} \, \varrho \, \int_Q [z_\varrho]^2 \, dx \, dt.
\end{equation}
The minimizer $z_\varrho \in L^2(Q)$ of \eqref{reduced L2 standard} is then
given as solution of the gradient equation
\begin{equation}\label{gradient equation L2 standard}
p_\varrho + \varrho \, z_\varrho = 0,
\end{equation}
where the adjoint $p_\varrho = {\mathcal{S}}^*(u_\varrho-\overline{u})$
solves the backward heat equation
\begin{equation}\label{backward heat equation}
  - \partial_t p_\varrho - \Delta_x p_\varrho = u_\varrho - \overline{u} \;
  \mbox{in} \; Q, \quad p_\varrho=0 \; \mbox{on} \; \Sigma, \quad
  p_\varrho=0 \; \mbox{on} \; \Sigma_T := \Omega \times \{ T \} .
\end{equation}
The resulting optimality system therefore consists of the primal
problem \eqref{heat equation}, the adjoint problem
\eqref{backward heat equation}, and the gradient equation
\eqref{gradient equation L2 standard}. While we can use this system
to determine the state $u_\varrho$, the control $z_\varrho$, and the
adjoint $p_\varrho$ simultaneously, one option is to eliminate the control
$z_\varrho$ from the gradient equation \eqref{gradient equation L2 standard}
first, and to consider a reduced optimality system to find
$(u_\varrho,p_\varrho)$ satisfying \eqref{backward heat equation} and
\begin{equation}\label{heat equation L2 standard}
  \partial_t u_\varrho - \Delta u_\varrho =
  - \frac{1}{\varrho} \, p_\varrho \; \mbox{in} \; Q, \quad
  u_\varrho = 0 \; \mbox{on} \; \Sigma, \quad u_\varrho = 0 \;
  \mbox{on} \; \Sigma_0 .
\end{equation}
When the adjoint $p_\varrho$ is known, the control $z_\varrho$ can be
computed by a simple post processing. As in \cite{LSTY_SISC:2021} we first
consider the simultaneous solution of the forward heat equation
\eqref{heat equation L2 standard} and the backward heat equation
\eqref{backward heat equation}, i.e., we consider a variational
formulation to find $u_\varrho \in X_0$ and
$p_\varrho \in X_T := \{ p \in Y : \partial_t \in Y^*, p=0 \; \mbox{on} \;
\Sigma_T \}$ such that
\begin{equation}\label{VF L2 standard}
  \begin{array}{rcl} \displaystyle
    \int_Q [\partial_t u_\varrho \, v +
    \nabla_x u_\varrho \cdot \nabla_x v ] \, dx \, dx
    & = & \displaystyle - \frac{1}{\varrho} \,
          \int_Q p_\varrho \, v \, dx \, dt
    \quad \mbox{for all} \; v \in Y, \\[4mm] \displaystyle
    \int_Q [- \partial_t p_\varrho \, q +
    \nabla_x p_\varrho \cdot \nabla_x q ] \, dx \, dt
    & = & \displaystyle
          \int_Q [u_\varrho - \overline{u}] \, q \, dx \, dt
          \quad \mbox{for all} \; q \in Y .
  \end{array}
\end{equation}
For the numerical solution of 
\eqref{VF L2 standard} we can use conforming space-time finite
element spaces $X_{0,h} = Y_{0,h} = S_h^1(Q) \cap X_0$ and
$X_{T,h} = Y_{T,h} = S_h^1(Q) \cap X_T$ of, e.g., piecewise linear,
continuous basis functions which are zero at the initial time $t=0$,
and at the final time $t=T$, respectively, and which are defined
with respect to some admissible decomposition of $Q$ into
shape-regular simplicial space-time finite elements.
Then the space-time finite element discretization of \eqref{VF L2 standard} is
to find $(u_{\varrho h}, p_{\varrho h}) \in X_{0,h} \times X_{T,h}$ such that
\begin{equation}\label{VF L2 standard FEM}
  \begin{array}{rcl} \displaystyle
    \int_Q [\partial_t u_{\varrho h} \, v_h +
    \nabla_x u_{\varrho h} \cdot \nabla_x v_h ] \, dx \, dx 
    & = & \displaystyle - \frac{1}{\varrho} \int_Q p_{\varrho h} \, v_h \,
          dx \, dt \quad \mbox{for all} \; v \in Y_{0,h}, \\[4mm]
    \displaystyle
    \int_Q [- \partial_t p_{\varrho h} \, q_h +
    \nabla_x p_{\varrho h} \cdot \nabla_x q_h ] \, dx \, dt
    & = & \displaystyle
          \int_Q [ u_{\varrho h} - \overline{u}] \, q_h \, dx \, dt
          \quad \mbox{for all} \; q_h \in Y_{T,h} .
  \end{array}
\end{equation}
The stability
and error analysis of this approach is given in \cite{LSTY_SISC:2021},
for general $\varrho \in {\mathbb{R}}_+$, and including a semi-linear
state equation with box constraints on the control.
The space-time finite element variational formulation
\eqref{VF L2 standard FEM} is equivalent to a linear system
of algebraic equations,
\begin{equation}\label{L2 standard System LGS}
  B_{0,h} \underline{u} = - \frac{1}{\varrho} \, M_h^\top \underline{p},
  \quad B_{T,h} \underline{p} = M_h \underline{u} - \underline{f},
\end{equation}
where all matrices and the vector $\underline{f}$ are defined according
to \eqref{VF L2 standard FEM}. In particular, the matrices
$B_{0,h}$ and $B_{T,h}$ representing the forward and backward
heat equation, respectively, are invertible.
Hence we can eliminate the discrete adjoint $\underline{p}$ to end
up with a Schur complement system to be solved,
\begin{equation}\label{L2 standard Schur LGS}
  \Big[ \varrho B_{0,h} +
  M_h^\top B_{T,h}^{-1} M_h \Big] \underline{u} =
  M_h^\top B_{T,h}^{-1} \underline{f} .
\end{equation}
Note that the mass matrix $M_h$ is in general not invertible; e.g., when the
numbers of degrees of freedom at the initial time $t=0$ and at the
final time $t=T$ are different, $M_h$ is rectangular.
While we have unique solvability of both the Schur complement system
\eqref{L2 standard Schur LGS}, and of the equivalent saddle point
system \eqref{L2 standard System LGS}, the construction of efficient
preconditioned iterative solution methods is at least not obvious.
The Schur complement matrix in \eqref{L2 standard Schur LGS} is
not symmetric, as one may expect when considering the minimization
of the reduced cost functional \eqref{reduced L2 standard}. This is
due to the fact that we have considered the discretization of the
forward and backward heat equations, but not of the optimality system.
Hence we lost the information that the heat equation was originally
a constraint for the minimization of \eqref{functional}. Instead of
\eqref{VF L2 standard} we now consider a variational formulation to
find $(u_\varrho,p_\varrho) \in X_0 \times Y$ such that
\begin{equation}\label{VF L2 richtig}
  \begin{array}{rcl}
    \displaystyle
    \frac{1}{\varrho} \, \int_Q p_\varrho \, q \, dx +
    \int_Q [\partial_t u_\varrho \, q + \nabla_x u_\varrho \cdot
    \nabla_x q] \, dx \, dt
    & = & 0 \quad \mbox{for all} \; q \in Y, \\[4mm]
    \displaystyle
    - \int_Q [p_\varrho \, \partial_t v + \nabla_x p_\varrho \cdot
    \nabla_x v ] \, dx \, dt +
    \int_Q u_\varrho \, v \, dx \, dt
    & = & \displaystyle \int_Q \overline{u} \, v \, dx \, dt \quad
          \mbox{for all} \; v \in X_0,
  \end{array}
\end{equation}
where we have used integration by parts in time to conclude the second
equation. For the numerical solution of \eqref{VF L2 richtig} we can use
again $X_{0,h} = S_h^1(Q) \cap X_0$, but now we consider the test space
$Y_h = S_h^1(Q) \cap Y$ without any conditions at the initial or final
time. Then the space-time finite element discretization of
\eqref{VF L2 richtig} is to find $(u_{\varrho h},p_{\varrho h}) \in
X_{0,h} \times Y_h$ such that
\begin{equation}\label{VF L2 richtig FEM}
  \begin{array}{rcl}
    \displaystyle
    \frac{1}{\varrho} \, \int_Q p_{\varrho h} \, q_h \, dx +
    \int_Q [\partial_t u_{\varrho h} \, q_h + \nabla_x u_{\varrho h} \cdot
    \nabla_x q_h] \, dx \, dt
    & = & 0, \\[4mm]
    \displaystyle
    - \int_Q [p_{\varrho h} \, \partial_t v_h + \nabla_x p_{\varrho h} \cdot
    \nabla_x v_h ] \, dx \, dt +
    \int_Q u_{\varrho h} \, v_h \, dx \, dt
    & = & \displaystyle \int_Q \overline{u} \, v_h \, dx \, dt 
  \end{array}
\end{equation}
is satisfied for all $(v_h,q_h) \in X_{0,h}\times Y_h$. This variational
formulation is equivalent to a linear system of algebraic equations,
\begin{equation}\label{L2 richtig System LGS}
  \frac{1}{\varrho} \, M_h \underline{p} + B_h \underline{u} =
  \underline{0}, \quad
  - B_h^\top \underline{p} + \overline{M}_h \underline{u} =
  \underline{f}
\end{equation}
where now the mass matrics $M_h$ and $\overline{M}_h$ are invertible,
but $B_h$ is rectangular and therefore not invertible. However, when
eliminating the discrete adjoint $\underline{p}$ this gives the
Schur complement system
\begin{equation}\label{L2 richtig Schur LGS}
  \Big[ \varrho \, B_h^\top M_h^{-1} B_h + \overline{M}_h \Big] \underline{u}
  = \underline{f},
\end{equation}
with a symmetric and positive definite Schur complement matrix
$S_h = \varrho \, B_h^\top M_h^{-1} B_h + \overline{M}_h$. This approach can
be analyzed similar as in \cite{LSY:2022}. While we can use a conjugate
gradient scheme to solve \eqref{L2 richtig Schur LGS}, we need to construct
a preconditioner for $S_h$ which is nothing than a space-time
finite element approximation of the bi-heat operator
$(-\partial_t - \Delta_x)(\partial_t - \Delta_x)$. Moreover, the application
of $S_h$ involves the inversion of the mass matrix, which can be replaced
by a lumped mass matrix as considered in \cite{LLSY_NLA:2023}
in the case of the Poisson equation.

The second approach becomes more intuitive when considering the
minimization of \eqref{functional} subject to the heat equation
\eqref{heat equation} in a different way. Instead of using the
solution operator $u_\varrho = {\mathcal{S}} z_\varrho$, we can write
the Dirichlet boundary value problem \eqref{heat equation} for the
heat equation as operator equation $B u_\varrho = z_\varrho$ where
we assume that $B:= \partial_t - \Delta_x : X \to L^2(Q)$ defines
an isomorphism, i.e.,
$X := \{ u \in X_0 : \partial_tu - \Delta_x u \in L^2(Q)\}$.
Following the abstract approach as considered in
\cite{LSY:2022} we introduce the reduced cost functional
\[
  \widehat{\mathcal{J}}(u_\varrho) \, = \,
  \frac{1}{2} \, \int_Q [u_\varrho - \overline{u}]^2 \, dx \, dt +
  \frac{1}{2} \, \varrho \, \| B u_\varrho \|^2_{L^2(Q)},
\]
whose minimizer $u_\varrho \in X$ is given as the unique solution
of the gradient equation
\begin{equation}\label{gradient L2 richtig}
  u_\varrho + \varrho B^* B u_\varrho = \overline{u} \, .
\end{equation}
When using the transformation $p_\varrho = - \varrho B u_\varrho$ we have
to solve the coupled system
\[
  B^* p_\varrho = u_\varrho - \overline{u}, \quad
  B u_\varrho = - \frac{1}{\varrho} \, p_\varrho,
\]
which is nothing than the coupled system
\eqref{backward heat equation} and \eqref{heat equation L2 standard}.
This shows the equivalence of both approaches, either using the
control-to-state map to eliminate the state, or to use the
partial differential equation to eliminate the control. Note that the
Schur complement system \eqref{L2 richtig Schur LGS} is a space-time
finite element approximation of the gradient equation
\eqref{gradient L2 richtig}, but since
$X_{0,h} = S_h^1(Q) \cap X_0 \not\subset X = \{ u \in X_0 :
\partial_t u - \Delta_x u \in L^2(Q) \}$, this is a non-conforming approach.
For a conforming discretization in the case of the Poisson equation, see
\cite{Brenner:2023}.

Since $B := \partial_t - \Delta_x
: X_0 \to Y^* = [L^2(0,T;H^1_0(\Omega))]^*$ defines an
isomorphism, we can also consider the control space
$Z = Y^*$. To realize the
norm $\| z \|_{Y^*} = \| w \|_Y = \| \nabla_x w \|_{L^2(Q)}$
we introduce $ w \in Y$ as the unique solution of the variational
formulation
\begin{equation}\label{Def A}
  \langle A w , v \rangle_Q :=
  \langle \nabla_x w , \nabla_x v \rangle_{L^2(Q)} =
  \langle z , v \rangle_Q \quad \mbox{for all} \; v \in Y.
\end{equation}
Then we can define the reduced cost functional
\[
  \overline{\mathcal{J}}(u_\varrho) =
  \frac{1}{2} \int_Q [u_\varrho - \overline{u}]^2 \, dx \, dt +
  \frac{1}{2} \, \varrho \, \langle B^* A^{-1} B u_\varrho, u_\varrho \rangle_Q,
\]
and its minimizer is given as the unique solution of the gradient equation
\begin{equation}\label{gradient equation Y}
  u_\varrho + \varrho B^* A^{-1} B u_\varrho = \overline{u} ,
\end{equation}
which we can write as
\begin{equation}\label{reduced optimality system Y}
  B^* p_\varrho = u_\varrho - \overline{u}, \quad
  \frac{1}{\varrho} \, A p_\varrho + B u_\varrho = 0 .
\end{equation}
This approach was first considered in \cite{GunzburgerKunoth:2011} using
wavelets for the discretization,
see also \cite{LSY:2022}. The space-time finite element discretization of
\eqref{reduced optimality system Y} is to find
$(u_{\varrho h},p_{\varrho h}) \in X_{0,h} \times Y_h$ such that
\begin{equation}\label{VF Bochner FEM}
  \begin{array}{rcl}
    \displaystyle
    \frac{1}{\varrho} \, \int_Q \nabla_x p_{\varrho h} \cdot
    \nabla_x q_h \, dx +
    \int_Q [\partial_t u_{\varrho h} \, q_h + \nabla_x u_{\varrho h} \cdot
    \nabla_x q_h] \, dx \, dt
    & = & 0, \\[4mm]
    \displaystyle
    - \int_Q [p_{\varrho h} \, \partial_t v_h + \nabla_x p_{\varrho h} \cdot
    \nabla_x v_h ] \, dx \, dt +
    \int_Q u_{\varrho h} \, v_h \, dx \, dt
    & = & \displaystyle \int_Q \overline{u} \, v_h \, dx \, dt 
  \end{array}
\end{equation}
is satisfied for all $(v_h,q_h) \in X_{0,h} \times Y_h$, and which is
equivalent to the coupled linear system
\begin{equation}\label{LGS System Y}
  \frac{1}{\varrho} \, A_h \underline{p} +
  B_h \underline{u} = \underline{0}, \quad
  \overline{M}_h \underline{u} - B_h^\top \underline{p} =
  \underline{f},
\end{equation}
where, in contrast to \eqref{L2 richtig System LGS}, the mass matrix $M_h$ is
replaced by the space-time finite element stiffness matrix of the spatial
Laplacian $A = - \Delta_x$. When eliminating $\underline{p}$ this now
results in the Schur complement system
\begin{equation}\label{LGS Schur Y}
  \Big[ \overline{M}_h + \varrho B_h^\top A_h^{-1} B_h \Big]
  \underline{u} = \underline{f} \, .
\end{equation}
The Schur complement matrix $S_h = \overline{M}_h + \varrho B_h^\top
A_h^{-1} B_h$ is again symmetric and positive definite, and we need to
construct a suitable preconditioner within the conjuate gradient
scheme. Due to the inversion of $A_h$ the application of
$S_h$ is now more involved. For a complete stability and error analysis
of this approach, see \cite{LSY:2022}. Note that we can also consider
the solution operator ${\mathcal{S}} : Y^* \to X_0$ in order to proceed as
in the first approach, i.e., similar as for
$\mathcal{S} : L^2(Q) \to X_0$. The
space-time finite element discretization of the forward and
backward heat equations within the related reduced optimality system
is given in \cite{LSTY_SINUM:2021}, but the final discrete system
shows less structure, similar as in \eqref{L2 standard Schur LGS}.

Note that
$S_h = B_h^\top A_h^{-1} B_h$ is a space-time finite element approximation
of the continuous Schur complement operator
$ S := B^* A^{-1} B : X_0 \to X_0^*$. 
The latter implies a norm in $X_0$, i.e., we have
\[
  \| u_\varrho \|^2_S = \langle A^{-1} B u_\varrho , B u_\varrho \rangle_Q =
  \langle A^{-1} z_\varrho , z_\varrho \rangle_Q = \| z_\varrho \|_{Y^*}^2 .
\]
Since the realization of $S$ includes the application of the inverse
$A^{-1}$, one may consider using a different operator implying a norm
in $X_0$ which is simpler
to handle, i.e., which avoids the use of any inverse operator.
Unfortunately, it is not obvious how to find a direct representation of
an equivalent norm in $X_0$.
However, the picture changes a lot when considering the variational
formulation \eqref{VF Bochner} in anisotropic Sobolev spaces.

\section{Optimal control in anisotropic Sobolev spaces}
\label{Section:Aniso}
Instead of using Bochner spaces within the variational formulation
\eqref{VF Bochner}, we now consider anisotropic Sobolev spaces
\cite{LionsMagenes:1972}. As in \cite{SteinbachZank:2020} we consider
a variational formulation to find $u_\varrho
\in H^{1,1/2}_{0;0,}(Q) := L^2(0,T;H^1_0(\Omega)) \cap
H^{1/2}_{0,}(0,T;L^2(\Omega))$ such that
\begin{equation}\label{VF Sobolev}
  \langle B u_\varrho , v \rangle_Q :=
  \langle \partial_t u_\varrho , v \rangle_Q +
  \langle \nabla_x u_\varrho , \nabla_x v \rangle_{L^2(Q)} =
  \langle z_\varrho , v \rangle_Q
\end{equation}
is satisfied for all $v \in H^{1,1/2}_{0;,0}(Q) :=
L^2(0,T;H^1_0(\Omega)) \cap H^{1/2}_{,0}(0,T;L^2(\Omega))$.
For $z_\varrho \in [H^{1,1/2}_{0;,0}(Q)]^*$ there exists a
unique solution $u_\varrho \in H^{1,1/2}_{0;0,}(Q)$ of the
variational formulation \eqref{VF Sobolev}, i.e., 
$B : H^{1,1/2}_{0;0,}(Q) \to [H^{1,1/2}_{0;,0}(Q)]^*$ is bijective.
We now define $w \in H^{1,1/2}_{0;,0}(Q)$ as unique
solution of the variational formulation
\begin{equation}\label{Def A aniso}
  \langle A w , v \rangle_Q :=
  \langle w , v \rangle_{H^{1,1/2}_{0;,0}(Q)} =
  \langle z , v \rangle_Q \quad \mbox{for all} \;
  v \in H^{1,1/2}_{0;,0}(Q).
\end{equation}
Although we are able to describe a computable representation of the
inner product in $H^{1,1/2}_{0;,0}(Q)$, this is not needed at this time.
Now we can choose $z_\varrho = Bu_\varrho \in Z := [H^{1,1/2}_{0;,0}(Q)]^*$
in order to conclude the reduced functional
\begin{equation}\label{reduced functional aniso inv}
  \check{\mathcal{J}}(u_\varrho) = \frac{1}{2} \int_Q
  [u_\varrho - \overline{u}]^2 \, dx \, dt +
  \frac{1}{2} \, \varrho \, \langle B^* A^{-1} B u_\varrho , u_\varrho \rangle_Q,
\end{equation}
whose minimizer formally satisfies the operator
equation \eqref{gradient equation Y}, and we may
proceed as above. But now the operator $S := B^* A^{-1} B :
H^{1,1/2}_{0;0,}(Q) \to [H^{1,1/2}_{0;0,}(Q)]^*$ implies a norm in
$H^{1,1/2}_{0;0,}(Q)$ for which we can find an equivalent norm which is
directly computable. In fact, following
\cite{SteinbachZank:2020}, a norm in $H^{1,1/2}_{0;0,}(Q)$ is given by
\begin{equation}\label{Def Norm}
  \| u \|_{H^{1,1/2}_{0;0,}(Q)}^2 :=
    \langle \partial_t u , {\mathcal{H}}_T u \rangle_Q +
    \| \nabla_x u \|^2_{L^2(Q)} =: \langle D u , u \rangle_Q
    = \| u \|_D^2 \, .
\end{equation}
Here we make use of the modified Hilbert transformation
${\mathcal{H}}_T$ as defined in \cite{SteinbachZank:2020}. For given
$u \in L^2(0,T)$ we consider the Fourier series
\[
  u(t) = \sum\limits_{k=0}^\infty u_k \sin \left( \Big(
    \frac{\pi}{2} + k \pi \Big) \frac{t}{T} \right), \quad
  u_k = \frac{2}{T} \int_0^T u(t)  \sin \left( \Big(
    \frac{\pi}{2} + k \pi \Big) \frac{t}{T} \right) \, dt,
\]
and define
\[
  {\mathcal{H}}_Tu(t) := \sum\limits_{k=0}^\infty u_k \cos \left( \Big(
    \frac{\pi}{2} + k \pi \Big) \frac{t}{T} \right) .
\]
Hence, instead of \eqref{reduced functional aniso inv} we finally
consider the reduced functional
\begin{equation}\label{reduced function aniso}
  \breve{\mathcal{J}}(u_\varrho) = \frac{1}{2} \int_Q
  [u_\varrho - \overline{u}]^2 \, dx \, dt +
  \frac{1}{2} \, \varrho \, \langle D u_\varrho , u_\varrho \rangle_Q .
\end{equation}
Note that
\[
  \| u_\varrho \|^2_D = \langle D u_\varrho , u_\varrho \rangle_Q =
  \langle D B^{-1} z_\varrho , B^{-1} z_\varrho \rangle_Q =
  \| z_\varrho \|^2_{B^{-1,*}D B^{-1}}
\]
defines an equivalent norm in $[H^{1,1/2}_{0;,0}(Q)]^*$. The minimizer of
\eqref{reduced function aniso} is now determined as unique solution
of the gradient equation
\begin{equation}\label{operator aniso}
u_\varrho + \varrho D u_\varrho = \overline{u} ,
\end{equation}
i.e., $u_\varrho \in H^{1,1/2}_{0;0,}(Q)$ solves the variational problem
\begin{equation}\label{VF aniso}
  \langle u_\varrho , v \rangle_{L^2(Q)} + \varrho \,
  \langle D u_\varrho , v \rangle_Q =
  \langle \overline{u} , v \rangle_{L^2(Q)} \quad
  \mbox{for all} \; v \in H^{1,1/2}_{0;0,}(Q) .
\end{equation}
Since this variational formulation corresponds to the abstract
formulation (2.9) in \cite{LSY:2022}, all regularization error estimates
as given in \cite[Lemma 3]{LSY:2022} remain valid.

\begin{lemma}
  Let $u_\varrho \in H^{1,1/2}_{0;0,}(Q)$ be the unique solution of the
  variational formulation \eqref{VF aniso}. For $\overline{u} \in L^2(Q)$
  there holds
  \begin{equation}\label{regularization error 0 0}
    \| u_\varrho - \overline{u} \|_{L^2(Q)} \leq \| \overline{u} \|_{L^2(Q)},
  \end{equation}
  while for $\overline{u} \in H^{1,1/2}_{0;0,}(Q)$ the following error
  estimates hold true:
  \begin{equation}
    \| u_\varrho - \overline{u} \|_{L^2(Q)} \leq \varrho^{1/2} \,
    \| \overline{u} \|_D,
  \end{equation}
  \begin{equation}
    \| u_\varrho - \overline{u} \|_D \leq \| \overline{u} \|_D.
  \end{equation}
  If in addition $D \overline{u} \in L^2(Q)$ is satisfied for
  $\overline{u} \in H^{1,1/2}_{0;0,}(Q)$,
  \begin{equation}\label{Regularisierung 0 2}
    \| u_\varrho - \overline{u} \|_{L^2(Q)} \leq \varrho \,
    \| D \overline{u} \|_{L^2(Q)}
  \end{equation}
  as well as
  \begin{equation}\label{Regularisierung 1 2}
    \| u_\varrho - \overline{u} \|_D \leq \varrho^{1/2} \,
    \| D \overline{u} \|_{L^2(Q)}
  \end{equation}
  follow. Moreover, when choosing $v = u_\varrho$ in
  \eqref{VF aniso}, this gives
\begin{equation}\label{regularization error 1 0}
  \| u_\varrho \|^2_{L^2(Q)} + \varrho \, \| u_\varrho \|_D^2 =
  \langle \overline{u}, u_\varrho \rangle_{L^2(Q)} \leq
  \| \overline{u} \|_{L^2(Q)} \| u_\varrho \|_{L^2(Q)} \leq
  \| \overline{u} \|^2_{L^2(Q)} .
\end{equation}

\end{lemma}

\begin{remark}
  From the definition \eqref{Def Norm} of $D$ and using the properties
  of the modified Hilbert transformation we conclude
  $D = {\mathcal{H}}_T^{-1} \partial_t - \Delta_x$, i.e.,
  \[
    \| D \overline{u} \|_{L^2(Q)} \leq
    \| {\mathcal{H}}_T^{-1} \partial_t \overline{u} \|_{L^2(Q)}
    + \| \Delta_x \overline{u} \|_{L^2(Q)} \leq
    \| \overline{u} \|_{H^1(0,T;L^2(\Omega))} +
    \| \overline{u} \|_{L^2(0,T;H^2(\Omega))} . 
  \]
  Hence it is sufficient to assume
  $\overline{u} \in H^{1,1/2}_{0;0,}(Q) \cap H^{2,1}(Q)$ to ensure
  \eqref{Regularisierung 0 2} and \eqref{Regularisierung 1 2},
  respectively.
\end{remark}

\begin{corollary}
  From $D = {\mathcal{H}}_T^{-1} \partial_t - \Delta_x :
  H^{1,1/2}_{0;0,}(Q) \cap H^{2,1}(Q) \to L^2(Q)$, and using
  \eqref{operator aniso}, we immediately have
$u_\varrho \in H^{2,1}(Q)$ for $\overline{u} \in L^2(Q)$,
and hence $z_\varrho = \partial_t u_\varrho - \Delta_x u_\varrho \in L^2(Q)$
follows.
\end{corollary}

\noindent
For the space-time finite element discretization of the variational
formulation \eqref{VF aniso}, let $X_h \subset H^{1,1/2}_{0;0,}(Q)$
be any conforming finite element space. The Galerkin space-time
finite element approximation of \eqref{VF aniso} is to find
$u_{\varrho,h} \in X_h$ such that
\begin{equation}\label{FEM aniso}
  \langle u_{\varrho,h} , v_h \rangle_{L^2(Q)} + \varrho \,
  \langle D u_{\varrho,h} , v_h \rangle_Q =
  \langle \overline{u} , v_h \rangle_{L^2(Q)} \quad
  \mbox{for all} \; v_h \in X_h .
\end{equation}
Using standard arguments we conclude Cea's lemma,
\begin{equation}\label{Cea}
  \| u_\varrho - u_{\varrho,h} \|^2_{L^2(Q)} + \varrho \,
  \| u_\varrho - u_{\varrho,h} \|^2_D \leq
  \| u_\varrho - v_h \|^2_{L^2(Q)} + \varrho \,
  \| u_\varrho - v_h \|^2_D \quad \mbox{for all} \; v_h \in X_h .
\end{equation}
In particular for $v_h \equiv 0$ we then obtain, using
\eqref{regularization error 1 0},
\[
  \| u_\varrho - u_{\varrho,h} \|^2_{L^2(Q)} \leq
  \| u_\varrho \|^2_{L^2(Q)} + \varrho \, \| u_\varrho \|^2_D \leq
  \| \overline{u} \|^2_{L^2(Q)},
\]
and hence, recall \eqref{regularization error 0 0},
\begin{equation}\label{FEM error 0 0}
  \| u_{\varrho,h} - \overline{u} \|_{L^2(Q)} \leq
  \| u_{\varrho,h} - u_\varrho \|_{L^2(Q)} +
  \| u_\varrho - \overline{u} \|_{L^2(Q)} \leq 2 \,
  \| \overline{u} \|_{L^2(Q)} .
\end{equation}
On the other hand, using the triangle inequality and
\eqref{Cea}, we also have
\begin{eqnarray*}
  \| u_{\varrho,h} - \overline{u} \|^2_{L^2(Q)}
  & \leq & 2 \, \| u_\varrho - \overline{u} \|^2_{L^2(Q)}
           + 2 \, \| u_\varrho - u_{\varrho,h} \|^2_{L^2(Q)} \\
  & \leq & 2 \, \| u_\varrho - \overline{u} \|^2_{L^2(Q)}
           + 2 \, \| u_\varrho - v_h \|^2_{L^2(Q)}
           + 2 \, \varrho \, \| u_\varrho - v_h \|_D^2 \\
  & \leq & 6 \, \| u_\varrho - \overline{u} \|_{L^2(Q)}^2
           + 4 \, \| \overline{u} - v_h \|^2_{L^2(Q)}
           + 4 \, \varrho \, \| u_\varrho - \overline{u} \|_D^2
           + 4 \, \varrho \, \| \overline{u} - v_h \|_D^2 .
\end{eqnarray*}
In particular when assuming
$\overline{u} \in H^{1,1/2}_{0;0,}(Q) \cap H^{2,1}(Q)$
we can use \eqref{Regularisierung 0 2} and
\eqref{Regularisierung 1 2} to conclude
\begin{equation}\label{Cea final}
  \| u_{\varrho,h} - \overline{u} \|^2_{L^2(Q)}
  \leq 10 \, \varrho^2 \, \| D \overline{u} \|_{L^2(Q)}^2
  + 4 \, \inf\limits_{v_h \in X_h}
  \Big[ \| \overline{u} - v_h \|^2_{L^2(Q)}
           + \varrho \, \| \overline{u} - v_h \|_D^2 \Big] .
\end{equation}
Hence it is sufficient to investigate the approximability of the
target $\overline{u}$ in the space-time finite element space $X_h$.
In particular when assuming
$\overline{u} \in H^{1,1/2}_{0;0,}(Q) \cap H^{2,1}(Q)$
this motivates the definition of tensor-product space-time finite
element spaces $X_h$ in order to derive approximation error estimates
which are anisotropic in space and time, see also
\cite{SteinbachZank:2020}.

Let $W_{h_x} = \mbox{span} \{ \psi_i \}_{i=1}^{M_x} \subset H^1_0(\Omega)$
be some spatial finite element space of piecewise linear basis functions
$\psi_i$ which are defined with respect to some admissible and globally
quasi-uniform finite element mesh with spatial mesh size $h_x$. Moreover,
$V_{h_t} := S_{h_t}^1(0,T) \cap H^{1/2}_{0,}(0,T) =
\mbox{span} \{ \varphi_k \}_{k=1}^{N_t}$ is the space of piecewise linear
functions, which are defined with respect to some uniform finite element
mesh with temporal mesh size $h_t$. Hence, we introduce the tensor-product 
space-time finite element space $X_h := W_{h_x} \otimes V_{h_t}$.

For a given $v \in H^{1/2}_{0,}(0,T;L^2(\Omega))$, we define the
$H^{1/2}_{0,}$ projection $Q_{h_t} v \in L^2(\Omega) \otimes V_{h_t}$
as the unique solution of the variational problem
\[
\langle \partial_t Q_{h_t}v , {\mathcal{H}}_T v_{h_t} \rangle_{L^2(Q)}
=
\langle \partial_t v , {\mathcal{H}}_T v_{h_t} \rangle_Q
\]
for all $v_{h_t} \in L^2(\Omega) \otimes V_{h_t}$. Moreover, for
$v \in L^2(0,T;H^1_0(\Omega))$, we define the $H^1_0$ projection
$Q_{h_x}v \in W_{h_x} \otimes L^2(0,T)$ as the unique solution of the
variational problem
\[
\int_0^T \int_\Omega \nabla_x Q_{h_x} v(x,t) \cdot
\nabla_x v_{h_x}(x,t) \, dx \, dt =
\int_0^T \int_\Omega \nabla_x v(x,t) \cdot
\nabla_x v_{h_x}(x,t) \, dx \, dt
\]
for all $v_{h_x} \in W_{h_x} \otimes L^2(0,T)$. It turns out that
$Q_{h_t}^{1/2} Q_{h_x}^1 v \in X_h$ is well--defined when
assuming $\partial_t v \in L^2(0,T;H^1_0(\Omega))$ and
$\nabla_x v \in H^{1/2}_{0,}(0,T;L^2(\Omega))$, and that the
projection operators $Q_{h_t}^{1/2}$, $Q_{h_x}^1$
and partial derivatives $\partial_t, \nabla_x$ commute in space and time,
see also \cite{SteinbachZank:2020,Zank:2019}.

\begin{lemma}\label{Lemma approximation aniso}
  Assume $u \in H^{1,1/2}_{0;0,}(Q) \cap H^{2,1}(Q)$, and $h_t \simeq h_x^2$.
  Then there hold the error estimates
  \begin{equation}\label{approximation 0 21}
    \| u - Q_{h_x} Q_{h_t} u \|_{L^2(Q)} \leq c \, h_x^2 \,
    | u |_{H^{2,1}(Q)},
  \end{equation}
  \begin{equation}\label{approximation 10 21}
    \| \nabla_x (u - Q_{h_x} Q_{h_t} u) \|_{L^2(Q)} \leq c \, h_x \,
    |u|_{H^{2,1}(Q)} ,
  \end{equation}
  and
  \begin{equation}\label{approximation D 21}
     \| u - Q_{h_x} Q_{h_t} u \|_{H^{1/2}_{0,}(0,T;L^2(\Omega))}
     \leq c \, h_x \, \Big[ \| u \|_{H^1(0,T;L^2(\Omega))} +
     \| u \|_{L^2(0,T;H^2(\Omega))} \Big] .
  \end{equation}
\end{lemma}

\proof{By the triangle inequality we first have, using standard
  stability and approximation error estimates for the $L^2$
  projections $Q_{h_x}$ and $Q_{h_t}$, respectively,
  \begin{eqnarray*}
    \| u - Q_{h_x} Q_{h_t} u \|_{L^2(Q)}
    & \leq & \| u - Q_{h_x} u \|_{L^2(Q)} +
             \| Q_{h_x} (u - Q_{h_t} u) \|_{L^2(Q)} \\
    & \leq & \| u - Q_{h_x} u \|_{L^2(Q)} +
             \| u - Q_{h_t} u \|_{L^2(Q)} \\
    & \leq & c_1 \, h_x^2 \, \| u \|_{L^2(0,T);H^2(\Omega))} +
             c_2 \, h_t \, \| u \|_{H^1(0,T;L^2(\Omega))} \\
    & \leq & c \, h_x^2 \, \Big[ \| u \|_{L^2(0,T);H^2(\Omega))} +
             \| u \|_{H^1(0,T;L^2(\Omega))} \Big] ,
  \end{eqnarray*}
  that is \eqref{approximation 0 21}. In the same way, but using an
  inverse inequality in $W_{h_x}$, we also obtain
  \begin{eqnarray*}
    \| \nabla_x (u - Q_{h_x} Q_{h_t} u) \|_{L^2(Q)}
    & \leq & \| \nabla_x (u - Q_{h_x} u) \|_{L^2(Q)} +
             \| \nabla_x Q_{h_x} (u - Q_{h_t} u) \|_{L^2(Q)} \\
    & \leq & \| \nabla_x (u - Q_{h_x} u) \|_{L^2(Q)} +
             c_{I} \, h_x^{-1} \, \| Q_{h_x} (u - Q_{h_t} u) \|_{L^2(Q)} \\
    & \leq & \| \nabla_x (u - Q_{h_x} u) \|_{L^2(Q)} +
             c_{I} \, h_x^{-1} \, \| u - Q_{h_t} u \|_{L^2(Q)} \\
    & \leq & c_1 \, h_x \, \| u \|_{L^2(0,T;H^2(\Omega))} +
             c_2 \, h_x^{-1} \, h_t \, \| u \|_{H^1(0,T;L^2(\Omega))} \\
    & \leq & c \, h_x \, \Big[ \| u \|_{L^2(0,T;H^2(\Omega))} +
             \| u \|_{H^1(0,T;L^2(\Omega))} \Big],
  \end{eqnarray*}
  that is \eqref{approximation 10 21}. We finally have, using an
  inverse inequality in $V_{h_t}$,
  \begin{eqnarray*}
    \| u - Q_{h_x} Q_{h_t} u \|_{H^{1/2}_{0,}(0,T;L^2(\Omega))} 
    & \leq & \| u - Q_{h_t} u \|_{H^{1/2}_{0,}(0,T;L^2(\Omega))} +
             \| Q_{h_t} (u - Q_{h_x} u) \|_{H^{1/2}_{0,}(0,T;L^2(Q))} \\
    & \leq & \| u - Q_{h_t} u \|_{H^{1/2}_{0,}(0,T;L^2(Q))} +
             c_{I} \, h_t^{-1/2} \, \| Q_{h_t}(u - Q_{h_x} u) \|_{L^2(Q)} \\
    & \leq & \| u - Q_{h_t} u \|_{H^{1/2}_{0,}(0,T;L^2(Q))} +
             c_{I} \, h_t^{-1/2} \, \| u - Q_{h_x} u \|_{L^2(Q)} \\
    & \leq & c_1 \, h_t^{1/2} \, \| u \|_{H^1(0,T;L^2(\Omega))} +
             c_2 \, h_t^{-1/2} \, h_x^2 \, \| u \|_{L^2(0,T;H^2(\Omega))} \\
    & \leq & c \, h_x \, \Big[ \| u \|_{H^1(0,T;L^2(\Omega))} +
             \| u \|_{L^2(0,T;H^2(\Omega))} \Big],
  \end{eqnarray*}
  i.e. \eqref{approximation D 21}.}

\noindent
Now, combining \eqref{Cea final} with the approximation error estimates
\eqref{approximation 0 21}, \eqref{approximation 10 21}, and
\eqref{approximation D 21}, this gives the following result.

\begin{theorem}\label{Thm error aniso}
  Let $u_{\varrho,h} \in X_h = W_{h_x} \otimes V_{h_t}$ be the unique solution
  of \eqref{FEM aniso}, where we assume $h_t \simeq h_x^2$. Assume
  $\overline{u} \in H^{1,1/2}_{0;0,}(Q) \cap H^{2,1}(Q)$, and consider
  $\varrho = h_x^2$. Then there holds
  \begin{equation}\label{final error 21}
    \| u_{\varrho,h} - \overline{u} \|_{L^2(Q)}
    \leq c \, h_x^2 \, \| \overline{u} \|_{H^{2,1}(Q)} .
\end{equation}
\end{theorem}

\noindent
As a consequence of the space-time finite element error estimates
\eqref{FEM error 0 0} and \eqref{final error 21}, and using a
space interpolation argument, we finally conclude the error estimate
\begin{equation}\label{final error aniso}
  \| u_{\varrho,h} - \overline{u} \|_{L^2(Q)}
  \leq c \, h_x^s \, \| \overline{u} \|_{H^{s,s/2}(Q)} 
\end{equation}
when assuming $\overline{u}\in H^{s,s/2}_{0;0,}(Q) =
[H^{1,1/2}_{0;0,}(Q),L^2(Q)]_s$ for $s\in[0,1]$, or
$\overline{u}\in H^{1,1/2}_{0;0,}(Q)\cap H^s(Q)$ for $s\in(1,2]$,
when using the parabolic scaling $h_t=h_x^2$, and $\varrho=h_x^2$.

The error estimate \eqref{final error 21} only assumes
$\overline{u} \in H^{1,1/2}_{0;0,}(Q) \cap H^{2,1}(Q)$, but requires
the parabolic scaling $h_t = h_x^2$. This is strongly related to
considering the control $z_\varrho \in L^2(Q)$ implying
$u_\varrho \in H^{1,1/2}_{0;0,}(Q) \cap H^{2,1}(Q)$ by maximal parabolic
regularity. But the error estimate \eqref{final error 21} was a
consequence of \eqref{Cea final}, where we have used the approximation
properties of the target $\overline{u}$ in the tensor-product
space-time finite element space $X_h = W_{h_x} \otimes V_{h_t}$.
Depending on the application in mind we may consider
$\overline{u}$ in anisotropic Sobolev spaces $H^{s,s/2}(Q)$, or
in isotropic spaces $H^s(Q)$. In particular, we now assume
$\overline{u} \in H^{1,1/2}_{0;0,}(Q) \cap H^2(Q)$. In this case we
can write the regularization error estimate
\eqref{Regularisierung 0 2} as
\begin{equation}\label{Regularisierung H2}
  \| u_\varrho - \overline{u} \|_{L^2(Q)} \leq \varrho \,
  \| \overline{u} \|_{H^2(Q)} .
\end{equation}
Moreover, we can adapt the space-time finite element error
estimates given in Lemma \ref{Lemma approximation aniso}
accordingly.

\begin{lemma}\label{Lemma approximation iso}
  Assume $u \in H^{1,1/2}_{0;0,}(Q) \cap H^2(Q)$, and $h_t \simeq h_x$.
  Then there hold the error estimates
  \begin{equation}\label{approximation 0 2}
    \| u - Q_{h_x} Q_{h_t} u \|_{L^2(Q)} \leq c \, h_x^2 \,
    | u |_{H^2(Q)},
  \end{equation}
  \begin{equation}\label{approximation 10 2}
    \| \nabla_x (u - Q_{h_x} Q_{h_t} u) \|_{L^2(Q)} \leq c \, h_x \,
    |u|_{H^2(Q)} ,
  \end{equation}
  and
  \begin{equation}\label{approximation D 2}
     \| u - Q_{h_x} Q_{h_t} u \|_{H^{1/2}_{0,}(0,T;L^2(\Omega))}
     \leq c \, h_x^{3/2} \, \| u \|_{H^2(Q)} .
  \end{equation}
\end{lemma}

\noindent
As a consequence, we have the following result:

\begin{corollary}
 Let $u_{\varrho,h} \in X_h = W_{h_x} \otimes V_{h_t}$ be the unique solution
  of \eqref{FEM aniso}, where we assume $h_t \simeq h_x$. Assume
  $\overline{u} \in H^{1,1/2}_{0;0,}(Q) \cap H^2(Q)$, and consider
  $\varrho = h_x^2$. Then there holds
  \begin{equation}\label{final error 2}
    \| u_{\varrho,h} - \overline{u} \|_{L^2(Q)}
    \leq c \, h_x^2 \, \| \overline{u} \|_{H^2(Q)} .
  \end{equation}
  Moreover, using a space interpolation argument, we also have
  \begin{equation}\label{final error s}
    \| u_{\varrho,h} - \overline{u} \|_{L^2(Q)}
    \leq c \, h_x^s \, \| \overline{u} \|_{H^s(Q)} ,
  \end{equation}
  when assuming $\overline{u}\in H^{s,s}_{0;0,}(Q) =
  [H^{1,1}_{0;0,}(Q),L^2(Q)]_s$ for $s\in[0,1]$,
  or $\overline{u}\in H^{1,1/2}_{0;0,}(Q)\cap H^s(Q)$ for $s\in(1,2]$.
\end{corollary}

\section{Optimal solution strategies}
\label{Section:Solver}
The Galerkin space-time variational formulation
\eqref{FEM aniso} is equivalent to a linear system of algebraic equations,
$K_h \underline{u} = \underline{f}$, where the symmetric and
positive definite stiffness matrix $K_h$
is given, due to the tensor-product ansatz space
$X_h = W_{h_x} \otimes V_{h_t}$, as
\begin{equation}\label{stiffness matrix}
  K_h = M_{h_t} \otimes M_{h_x} + \varrho \, \Big[
  A_{h_t} \otimes M_{h_x} + M_{h_t} \otimes A_{h_x} \Big] \in
  {\mathbb{R}}^{N_t \cdot M_x \times N_t \cdot M_x} ,
\end{equation}
where
\[
  A_{h_t}[j,i] = \langle \partial_t \varphi_i , {\mathcal{H}}_T
  \varphi_j \rangle_{L^2(0,T)}, \quad
  M_{h_t}[j,i] = \langle \varphi_i , \varphi_j \rangle_{L^2(0,T)}
  \quad \mbox{for} \; i,j=1,\ldots,N_t,
\]
and
\[
  A_{h_x}[\ell,k] = \langle \nabla_x \psi_k , \nabla_x \psi_\ell
  \rangle_{L^2(\Omega)}, \quad
  M_{h_x}[\ell,k] = \langle \psi_k , \psi_\ell \rangle_{L^2(\Omega)}
  \quad \mbox{for} \; k,\ell=1,\ldots,M_x.
\]

\begin{lemma}\label{lem:spectralEquivalence}
  For the finite element stiffness matrix $K_h$ as given in
  \eqref{stiffness matrix}, and for the mass matrix
  $M_h = M_{h_t} \otimes M_{h_x}$ there hold the spectral
  equivalence inequalities
  \begin{equation}
    ( M_h \underline{u} , \underline{u} ) \leq
    (K_h \underline{u} , \underline{u} ) \leq c \,
    ( M_h \underline{u} , \underline{u} ) \quad
    \mbox{for all} \; \underline{u} \in
    {\mathbb{R}}^{N_t \cdot M_x},
  \end{equation}
  when choosing $\varrho = h_x^2$ and $h_t\simeq h_x^2$ or $h_t\simeq h_x$.
\end{lemma}

\proof{While the lower estimate is trivial, to prove the upper one
  we use the finite element isomorphism $\underline{u} \in
  {\mathbb{R}}^{N_t \cdot M_x} \leftrightarrow u_h \in X_h =
  W_{h_x} \otimes V_{h_t}$ to write, using inverse inequalities in
  $V_{h_t}$ and $W_{h_x}$, respectively,
  \begin{eqnarray*}
    (K_h \underline{u} , \underline{u} )
    & = & \| u_h \|^2_{L^2(Q)} + \varrho \,
          \Big[ \| u_h \|^2_{H^{1/2}_{0,}(0,T;L^2(\Omega))} +
          \| \nabla_x u_h \|^2_{L^2(\Omega)} \Big] \\
    & = & \| u_h \|^2_{L^2(Q)} + h_x^2 \,
          \Big[ c_1 \, h_t^{-1} \| u_h \|^2_{L^2(Q)} + c_2 \, h_x^{-2}
          \| u_h \|^2_{L^2(\Omega)} \Big] \\
    & \leq & c \, (1+h_x^2 \, h_t^{-1}) \, \| u_h \|^2_{L^2(Q)} \, \leq \,
             c \, (M_h \underline{u},\underline{u}) \, .
  \end{eqnarray*}
  While the last estimate is obvious for $h_t \simeq h_x^2$, for
  $h_t \simeq h_x$ we use $h_x^2 h_t^{-1} \simeq h_x \leq c$.}

\noindent
Since in the case of globally quasi-uniform meshes, the mass matrix
$M_h = M_{h_t} \otimes M_{h_x}$ is spectrally equivalent to the identity,
we can use a conjugate gradient scheme without preconditioner to solve
the linear system $K_h \underline{u}=\underline{f}$.

However, due to the structure of the stiffness matrix $K_h$ as given
in \eqref{stiffness matrix}, and following
\cite{LangerZank:2021}, we are able to construct an efficient direct
solver. For this we need to have all eigenvalues and eigenvectors of
the general eigenvalue problem 
$A_{h_t}\underline{v}= \lambda M_{h_t}\underline{v}$, which were computed
numerically in \cite{LangerZank:2021}.
Following \cite{LSZ:2024}, we are now able to derive explicit
representations for all eigenvectors and eigenvalues of
$A_{h_t}$.

\begin{lemma}\label{lem:eigenvalues}
  The solution of the generalized eigenvalue problem
  \begin{equation}\label{gen EVP}
    A_{h_t} \underline{v} = \lambda \, M_{h_t} \underline{v}
  \end{equation}
  is given by the eigenvectors $\underline{v}_\ell$, $\ell=0,\ldots,N_t-1$,
  with components
  \begin{equation}\label{EV}
    v_{\ell,i} = \sin \left( \Big( \frac{\pi}{2} + \ell \pi \Big)
      \frac{i}{N_t} \right) \quad \mbox{for} \; i=1,\ldots,N_t,
  \end{equation}
  and with the corresponding eigenvalue
  \begin{equation}\label{Eigenwert}
    \lambda_\ell =
    \frac{3\pi}{2T} \frac{\sin^4 x_\ell}{x_\ell^4}
        \, \frac{1}{2+\cos \left( 2 x_\ell \right)} \,
        \sum\limits_{\mu=0}^\infty \left[
        \frac{(2\ell+1)^4}{(4\mu N_t + 2\ell+1)^3}
        +
        \frac{(2\ell+1)^4}{(4\mu N_t + 4N_t-1-2\ell)^3}\right] ,
    \end{equation}
    where
    \[
      x_\ell = \Big( \frac{\pi}{2} + \ell \pi \Big) \frac{1}{2N_t} .
  \]
\end{lemma}

\proof{For $\ell=0,\ldots,N_t-1$, let
  $\underline{v}_\ell = ( v_{\ell,i} )_{i=1}^{N_t}$
  be the eigenvector as given in \eqref{EV}, where in addition we can
  introduce $v_{\ell,0}=0$. When considering the
  matrix vector product $M_{h_t} \underline{v}_\ell$, we distinguish two cases:
  For $j=1,\ldots,N_t-1$ we have, recall that the temporal mesh is
  uniform,
  \begin{eqnarray*}
  (M_{h_t} \underline{v}_\ell)_j
    & = & \frac{1}{6} \, h_t \, \Big[ v_{\ell,j-1} + 4 v_{\ell,j} + v_{\ell,j+1}
          \Big] \\
  &&  \hspace*{-2cm} = \, \frac{1}{6} \, h_t \, \left[
        \sin \left(
        \Big( \frac{\pi}{2} + \ell \pi \Big) \frac{j-1}{N_t} \right)
        +4 \sin \left(
        \Big( \frac{\pi}{2} + \ell \pi \Big) \frac{j}{N_t} \right)
        + \sin \left(
        \Big( \frac{\pi}{2} + \ell \pi \Big) \frac{j+1}{N_t} \right)
        \right] \\
  && \hspace*{-2cm} = \, \frac{1}{6} \, h_t \, \left[
        4 \sin \left(
        \Big( \frac{\pi}{2} + \ell \pi \Big) \frac{j}{N_t} \right)
        + 2 \sin \left(
        \Big( \frac{\pi}{2} + \ell \pi \Big) \frac{j}{N_t} \right)
        \cos \left(
        \Big( \frac{\pi}{2} + \ell \pi \Big) \frac{1}{N_t} \right)
        \right] \\
  && \hspace*{-2cm} = \, \frac{1}{3} \, h_t \, \left[
        2 + \cos \left(
        \Big( \frac{\pi}{2} + \ell \pi \Big) \frac{1}{N_t} \right)
        \right]
        \sin \left(
     \Big( \frac{\pi}{2} + \ell \pi \Big) \frac{j}{N_t} \right) \\
    && \hspace*{-2cm} = \, \frac{1}{3} \, h_t \, \Big[ 2 + \cos 2x_\ell
       \Big] \, \sin \left(
     \Big( \frac{\pi}{2} + \ell \pi \Big) \frac{j}{N_t} \right),
  \end{eqnarray*}
  while for $j=N_t$ we have
  \begin{eqnarray*}
  (M_{h_t} \underline{v}_\ell)_{N_t}
  & = & \frac{1}{6} \, h_t \, \Big[ v_{\ell,N_t-1} + 2 v_{\ell,N_t} \Big] \\
  && \hspace*{-2cm} = \, \frac{1}{6} \, h_t \, \left[
        \sin \left(
        \Big( \frac{\pi}{2} + \ell \pi \Big) \frac{N_t-1}{N_t} \right)
        + 2 \sin 
        \Big( \frac{\pi}{2} + \ell \pi \Big) 
        \right] \\
  && \hspace*{-2cm} = \, \frac{1}{6} \, h_t \, \left[
        \sin 
        \Big( \frac{\pi}{2} + \ell \pi \Big) 
        \cos \left(
        \Big( \frac{\pi}{2} + \ell \pi \Big) \frac{1}{N_t} \right)
        -
        \cos 
        \Big( \frac{\pi}{2} + \ell \pi \Big) 
        \sin \left(
     \Big( \frac{\pi}{2} + \ell \pi \Big) \frac{1}{N_t} \right)
     \right. \\
    && \left. \hspace*{9.5cm} 
        + 2 \sin 
        \Big( \frac{\pi}{2} + \ell \pi \Big) 
        \right] \\
  && \hspace*{-2cm} = \, \frac{1}{6} \, h_t \, \left[
        \cos \left(
        \Big( \frac{\pi}{2} + \ell \pi \Big) \frac{1}{N_t} \right)
     + 2 \right]
     \sin \Big( \frac{\pi}{2} + \ell \pi \Big) \\
    && \hspace*{-2cm} = \, \frac{1}{6} \, h_t \, \Big[
       2 + \cos 2 x_\ell \Big]
     \sin \Big( \frac{\pi}{2} + \ell \pi \Big) \, .
\end{eqnarray*}
Now we are going to prove that $\underline{v}_\ell$ is indeed an
eigenvector of the generalized eigenvalue problem \eqref{gen EVP}.
For this we need to compute
\[
  (A_{h_t} \underline{v}_\ell)_j =
  \sum\limits_{i=1}^{N_t} A_{h_t}[j,i] v_{\ell,i} =
  \sum\limits_{i=1}^{N_t} v_{\ell,i} \, \langle \partial_t \varphi_i ,
  {\mathcal{H}}_T \varphi_j \rangle_{L^2(0,T)} =
  \langle \partial_t v_{\ell,h} , {\mathcal{H}}_T \varphi_j
  \rangle_{L^2(0,T)},
\]
where
\[
  v_{\ell,h}(t) = \sum\limits_{i=1}^{N_t} v_{\ell,i} \varphi_i(t) =
  \sum\limits_{k=0}^\infty A_k \sin \left( \Big(
    \frac{\pi}{2} + k\pi \Big) \frac{t}{T} \right) ,
\]
and with the Fourier coefficients
\begin{eqnarray*}
  A_k
  & = & \frac{2}{T} \int_0^T v_{\ell,h}(t) \,
        \sin \left( \Big( \frac{\pi}{2} + k\pi \Big) \frac{t}{T} \right) \\
  & = & \frac{2}{T} \sum\limits_{i=1}^{N_t} v_{\ell,i} \int_0^T \varphi_i(t) \,
        \sin \left( \Big( \frac{\pi}{2} + k\pi \Big) \frac{t}{T} \right)
        \, dt \, = \, \sum\limits_{i=1}^{N_t} v_{\ell,i} A_k^i \, .
\end{eqnarray*}
Note that all following more technical computations follow similar
as in \cite{LSZ:2024}, where we have considered the Fourier expansion
of piecewise constant basis functions. In particular, for $i=1,\ldots,N_t-1$,
\[
  A_k^i
  \, = \, \frac{2}{T} \int_0^T \varphi_i(t) \,
  \sin \left( \Big( \frac{\pi}{2} + k\pi \Big) \frac{t}{T} \right) \, dt
  \, = \, \frac{2}{N_t} \, \frac{\sin^2 x_k}{x_k^2} \,
  \sin \left( \Big( \frac{\pi}{2} + k \pi \Big) \frac{i}{N_t} \right),
\]
and
\[
  A_k^{N_t}
  \, = \, \frac{2}{T} \int_0^T \varphi_{N_t}(t) \,
  \sin \left( \Big( \frac{\pi}{2} + k\pi \Big) \frac{t}{T} \right) \, dt
  \, = \, \frac{1}{N_t} \, \frac{\sin^2 x_k}{x_k^2} \,
  \sin \left( \frac{\pi}{2} + k \pi \right) .
\]    
Hence we can write
\[
  \partial_t v_{\ell,h}(t) = \frac{1}{T} \sum\limits_{k=0}^\infty A_k
  \left( \frac{\pi}{2} + k\pi \right)
  \cos \left( \Big( \frac{\pi}{2} + k\pi \Big) \frac{t}{T} \right)
\]
and
\[
  {\mathcal{H}}_T \varphi_j(t) =
  \sum\limits_{\ell=0}^\infty A_\ell^j
  \cos \left( \Big( \frac{\pi}{2} + \ell\pi \Big) \frac{t}{T} \right)
\]
to conclude
\begin{eqnarray*}
  && \langle \partial_t v_{\ell,h} ,
     {\mathcal{H}}_T \varphi_j \rangle_{L^2(0,T)} \\
  && \hspace*{1cm} 
     = \, \frac{1}{T} \sum\limits_{k=0}^\infty \sum\limits_{\ell=0}^\infty
     A_k A_\ell^j \left( \frac{\pi}{2} + k\pi \right) \int_0^T
     \cos \left( \Big( \frac{\pi}{2} + k\pi \Big) \frac{t}{T} \right)
     \cos \left( \Big( \frac{\pi}{2} + \ell\pi \Big) \frac{t}{T} \right)
     \, dt \\
  && \hspace*{1cm} = \, \frac{1}{2} \sum\limits_{k=0}^\infty A_k A_k^j
     \left( \frac{\pi}{2} + k\pi \right) \, = \,
     \frac{1}{2} \sum\limits_{i=1}^{N_t} v_{\ell,i} \sum\limits_{k=0}^\infty
     A_k^i A_k^j \left( \frac{\pi}{2} + k\pi \right) \, .
\end{eqnarray*}
For all $j=1,\ldots,n$ and for all $k=0,\ldots,2n-1$ we obtain
the recurrence relation
\[
  A_{k+2\mu N_t}^j = \frac{(2k+1)^2}{(4\mu N_t+2k+1)^2} \, A_k^j \quad
  \mbox{for} \; \mu \in {\mathbb{N}},
\]
while for $k=0,\ldots,N_t-1$ we conclude
\[
A_{2N_t-1-k}^j = - \frac{(2k+1)^2}{(4N_t-1-2k)^2} \, A_k^j \, .
\]
Hence we can write
\[
  \sum\limits_{k=0}^\infty
  A_k^i A_k^j \left( \frac{\pi}{2} + k\pi \right) \, = \,
  \sum\limits_{k=0}^{N_t-1} \gamma(k,N_t) A_k^i A_k^j ,
\]
where
\[
\gamma(k,N_t) = \frac{\pi}{2} 
\sum\limits_{\mu=0}^\infty
     \left[
        \frac{(2k+1)^4}{(4\mu N_t + 2k+1)^3}
        +
     \frac{(2k+1)^4}{(4\mu N_t + 4n-1-2k)^3}\right] \, .
\]
When using
\[
  \sum\limits_{i=1}^{N_t-1}  \sin \left(
    \Big( \frac{\pi}{2} + \ell \pi \Big) \frac{i}{N_t} \right)
  \sin \left(
    \Big( \frac{\pi}{2} + k\pi \Big) \frac{i}{N_t} \right)
  +
  \frac{1}{2} \sin \left( \frac{\pi}{2} + \ell\pi \right)
  \sin \left( \frac{\pi}{2} + k\pi \right) =
  \frac{N_t}{2} \, \delta_{k\ell},
\]
we now conclude, for $j=1,\ldots,N_t$,
\begin{eqnarray*}
  (A_{h_t}\underline{v}_\ell)_j
  & = & \frac{1}{2} \sum\limits_{i=1}^{N_t} v_{\ell,i}
        \sum\limits_{k=0}^{N_t-1} \gamma(k,N_t) A_k^i A_k^j \\
  & = & \frac{1}{2} \sum\limits_{k=0}^{N_t-1} \gamma(k,N_t) \left[
        \sum\limits_{i=1}^{N_t-1} v_{\ell,i} A_k^i +
        v_{\ell,N_t} A_k^{N_t} \right] A_k^j \\
  & = & \frac{1}{N_t}
        \sum\limits_{k=0}^{N_t-1} \gamma(k,N_t)
        \frac{\sin^2x_k}{x_k^2} \left[
        \sum\limits_{i=1}^{N_t-1} 
        \sin \left( \Big( \frac{\pi}{2} + k\pi \Big) \frac{i}{N_t}
        \right)
        \sin \left( \Big( \frac{\pi}{2} + \ell \pi \Big) \frac{i}{N_t}
        \right) \right. \\
  && \left. \hspace*{6cm}
     + 
     \sin \left( \frac{\pi}{2} + k\pi \right)
        \sin \left( \Big( \frac{\pi}{2} + \ell \pi \Big) \frac{i}{N_t}
        \right)\right] A_k^j \\
  & = & \frac{1}{2} \gamma(\ell,N_t)
        \frac{\sin^2x_\ell}{x_\ell^2} A_\ell^j \, .
\end{eqnarray*}
In particular for $j=1,\ldots,N_t-1$ this gives
\[
  (A_{h_t}\underline{v}_\ell)_j
  = \frac{1}{N_t} \gamma(\ell,N_t) \frac{\sin^4x_\ell}{x_\ell^4}
        \sin \left( \Big( \frac{\pi}{2} + \ell \pi \Big) \frac{j}{N_t}
        \right) 
  = \frac{3}{T} \gamma(\ell,N_t) \frac{\sin^4x_\ell}{x_\ell^4}
        \frac{1}{2 + \cos 2 x_\ell} \, (M_{h_t} \underline{v}_\ell)_j,
\]
while for $j=N_t$ we have
\[
  (A_{h_t}\underline{v}_\ell)_{N_t} =
  \frac{1}{2N_t} \gamma(\ell,N_t) \frac{\sin^4x_\ell}{x_\ell^4}
  \sin \left( \frac{\pi}{2} + \ell \pi \right) =
  \frac{3}{T} \gamma(\ell,N_t) \frac{\sin^4x_\ell}{x_\ell^4}
  \frac{1}{2+\cos 2x_\ell} \,
  (M_{h_t} \underline{v}_\ell)_{N_t}.
\]
In both cases, this is  $ (A_{h_t}\underline{v}_\ell)_j =
\lambda_\ell (M_{h_t}\underline{v}_\ell)_j$,
$j=1,\ldots,N_t$, with $\lambda_\ell$ as given in
\eqref{Eigenwert}.}

\noindent
Now we are in the position to describe the solution of the linear
system $K_h \underline{u}=\underline{f}$ where the stiffness matrix
$K_h$ is given as in \eqref{stiffness matrix}. When we multiply the
system from the left by $M_{h_t}^{-1} \otimes I_{h_x}$ this gives
\begin{equation}
  \left(
    I_{h_t} \otimes M_{h_x} + \varrho \, M_{h_t}^{-1} A_{h_t} \otimes M_{h_x} +
    \varrho \, I_{h_t} \otimes A_{h_x}\right) \underline{u} =
  \left(
    M_{h_t}^{-1} \otimes I_{h_x} \right) \underline{f}.
  \label{eq:directsolver:firstStep}
\end{equation}
Let $C_{h_t}$ be the matrix containing the eigenvectors $\underline{v}_\ell$
of the matrix $M_{h_t}^{-1}A_{h_t}$ as columns, and $D_{h_t}$ be the diagonal
matrix containing the corresponding eigenvalues $\lambda_\ell$ as given in
Lemma \ref{lem:eigenvalues}. With this we can bring the system
\eqref{eq:directsolver:firstStep} into time-diagonal form by defining
\[
    \underline{u} := \left( C_{h_t} \otimes I_{h_x} \right) \underline{v},
\]
and by multiplying \eqref{eq:directsolver:firstStep} from the left by
$C_{h_t}^{-1} \otimes I_{h_x}$ to obtain
\[
  \left(
    I_{h_t} \otimes M_{h_x} + \varrho \, D_{h_t} \otimes M_{h_x} +
    \varrho \, I_{h_t} \otimes A_{h_x} \right) \underline{v} =
  \left( C_{h_t}^{-1} M_{h_t}^{-1} \otimes I_{h_x}\right)
  \underline{f} =: \underline{g}.
\]
Due to the space-time tensor-product structure, and introducing
$\underline{v}^{(i)}$ and $\underline{g}^{(i)}$ as the restriction
of $\underline{v}$ and $\underline{g}$ to the discrete time $t_i$, it
remains to solve $N_t$ independent systems
\begin{equation}
  \left(
    M_{h_x} + \varrho \, \lambda_i \, M_{h_x} + \varrho \, A_{h_x} \right)
  \underline{v}^{(i)} = \underline{g}^{(i)}, \quad i=1,\ldots,N_t.
  \label{eq:directsolver:spatialProblems}
\end{equation}
The system matrices of the spatial problems
\eqref{eq:directsolver:spatialProblems} are all
symmetric, positive definite, and spectrally equivalent to the mass matrix
$M_{h_x}$, which can be shown analogously to
Lemma \ref{lem:spectralEquivalence} when choosing $\varrho = h_x^2$. Since
we use a globally quasi-uniform spatial mesh, we can use a conjugate gradient
scheme without preconditioner for the iterative solution of
\eqref{eq:directsolver:spatialProblems}. Note that $M_{h_t}$ is of
tridiagonal form, and thus can be inverted efficiently with linear complexity.
So the only two performance bottlenecks are the inversion of $C_{h_t}$ as well
as the multiplication by $C_{h_t}$ as these are dense matrices. In summary,
at this time we end up with a complexity of $\mathcal{O}(N_t^2M_x)$.
However, when taking a closer look on the matrix vector product
\[
  w_i = \left( C_{h_t} \underline {v}\right)_i =
  \sum_{k=0}^{N_t-1} \sin \left(\left(\frac{\pi}{2}+k \pi\right)
    \frac{i}{N_t}\right) v_k, \quad i = 1, \ldots, N_t,
\]
we notice, that it is two times the type two discrete sine transform of
$\vec{v}$. The discrete sine transform as well as its inverse are closely
related to the discrete Fourier transform, which is well understood and can be computed
efficiently with $\mathcal{O}(N_t \log N_t)$ complexity \cite{FFTW05}. This is
the approach we use in our implementation. Hence we end up with an
optimal complexity of ${\mathcal{O}}(M_x N_t \log N_t)$, which is de facto
linear in the number of degrees of freedom.

\section{Numerical results}\label{Section:Numerics}
 \input{includes/numerical_examples.tex}

\section{Conclusions}
In this paper we have formulated and analyzed a new approach to solve
distributed optimal control problems subject to the heat equation as
a minimization problem with respect to the state in the anisotropic
Sobolev space $H^{1,1/2}_{0;0,}(Q)$ where we also present an efficient
approach of optimal complexity to evaluate the norm in
$H^{1,1/2}_{0;0,}(Q)$. While the application of the fast sine
transformation requires a space-time tensor product mesh which is uniform
in time, we can use either tensor-product or simplicial meshes in space.
In order to handle adaptive meshes in space, and as in 
\cite{LLSY:adaptive} for the elliptic case, we can introduce a variable
energy regularization in space, i.e., instead of
\eqref{reduced function aniso} we can minimize the reduced cost functional
\[
  \widetilde{\mathcal{J}}(u_\varrho) =
  \frac{1}{2} \, \| u_\varrho - \overline{u} \|^2_{L^2(Q)} +
  \frac{1}{2} \, \varrho_t \, \langle \partial_t u_\varrho ,
  {\mathcal{H}}_T u_\varrho \rangle_Q +
  \frac{1}{2} \, \int_0^T \int_\Omega \varrho_x(x) \,
  |\nabla_x u_\varrho(x,t)|^2 \, dx \, dt,
\]
where $\varrho_t$ is a constant to be chosen accordingly, and
$\varrho_x(x) = h_{x,\ell}^2$ for $x \in \tau_\ell \subset \Omega$.

The focus of this paper was on the efficient solution of the optimal
control problem without additional constraints neither on the state, nor
on the control. But we showed that state constaints can be incorporated, see also \cite{löscher2024optimalcomplexitysolutionspacetime}. The same is true for control constraints, see \cite{GanglLS:2023}. Both lead to variational
inequalities to be solved by some iterative procedure, which again
require an efficient handling of all matrices as described in this
paper. Related results of these extensions will be published elsewhere.

It is without saying, that the proposed approach is also of interest to
other optimal control and inverse problems subject to parabolic
evolution equations, including nonlinear problems.

\bigskip
\noindent
\textbf{Acknowledgement:} 
This work has been supported by the Austrian
Science Fund (FWF) under the Grant Collaborative Research Center
TRR361/F90: CREATOR Computational Electric Machine Laboratory.

\end{document}

%% file: includes/numerical_examples.tex
For the numerical experiments, let $\Omega = (0,1)^3$, and $T=1$, i.e.,
$Q = (0,1)^4$. For simplicity, we also use a uniform tensor product 
mesh with $n_x$ finite elements of mesh size $h_x=1/n_x$ in each
coordinate for the spatial resolution. For $T=1$, the time interval
$(0,1)$ is decomposed into $n_t$ temporal finite elements of mesh size
$h_t=1/n_t$. In the case of a globally uniform mesh we have $n_t=n_x=N_t$,
and therefore, due to homogeneous initial and boundary conditions,
$M_x\cdot N_t = (n_x-1)^3 n_x$ degrees of freedom. In the case of parabolic scaling we have
$n_t = n_x^2$, and therefore, $(n_x-1)^3 n_x^2$ degrees of freedom.
Note that instead of a spatial tensor product mesh we may also consider
a globally quasi-uniform simplicial mesh, but in time we have to use
a uniform mesh in order to characterize the eigenvalues and eigenvectors
of $A_{h_t}$ as given in Lemma \ref{lem:eigenvalues}. Otherwise, the
symmetric and positive definite but dense matrix $A_{h_t}$ has to be
factorized numerically.

While for the relaxation parameter $\varrho$ we conclude the optimal
choice $\varrho = h_x^2$ in all cases, the order of convergence depends on the
regularity of the target $\overline{u}$.

As a first example we consider a smooth target
$\overline{u}_s \in H^{1,1/2}_{0;0,}(Q) \cap H^2(Q)$ also satisfying
homogeneous initial and boundary conditions,
\[
  \overline{u}_s(x, t) = t^2 \, x_1 \, (1 - x_1) \,
  x_2 \, (1 - x_2) \, x_3 \, (1 - x_3 ) .
\]
The numerical results for this example are given in
Table \ref{tab:uS} where we observe a second order convergence (eoc),
as predicted. The computing time in ms covers the solution of all
involved linear systems, where the relative accuracy of the conjugate
gradient scheme was set to $10^{-12}$. The ratio of the number of
degrees of freedom over time is almost constant, indicating linear complexity. Recall that the
application of the fast sine transformation includes a factor which is
logarithmic in $n_t$. All computations were done
on a Intel Xeon E5-2630 v3 CPU (20~MB cache, 3.20~GH) with 256~GB DDR4 RAM on a single thread.

\begin{table}[H]
    \centering
    \begin{filecontents*}{simulation/h22_3d/transformed_data.csv}
dof,simulationTime,spatialProblemsTime,percentageSpatialSolve,nx,nt,eoc,l2Error,timePerDof
2,6,0,0.0,2,2,0.0,0.00246325,3000.0
108,18,7,38.888888888888886,4,4,0.420255,0.00184077,166.66666666666666
2744,151,112,74.17218543046357,8,8,1.00915,0.000914567,55.029154518950435
54000,1623,1312,80.83795440542205,16,16,1.58096,0.000305704,30.055555555555557
953312,25628,20908,81.58264398314344,32,32,1.86325,8.40247e-05,26.8831190628042
16003008,419179,346226,82.5962178448825,64,64,1.95822,2.16234e-05,26.193763072542364
262193024,6656635,5499839,82.62191031955335,128,128,1.98666,5.45608e-06,25.38829942325239
\end{filecontents*}
\pgfplotstabletypeset[%
    col sep=comma, %
    columns={dof,nx,nt,l2Error,eoc, simulationTime, timePerDof}, %
    every head row/.style={before row=\toprule,after row=\midrule}, %
      every last row/.style={after row=\bottomrule}, %
      columns/simulationTime/.style={column name={time [ms]}, column type={r}, fixed}, %
      columns/timePerDof/.style={column name={time/dof [ns]}, column type={r}, fixed}, %
      columns/dof/.style={column name={DoF}, column type={r}, fixed}, %
      columns/nx/.style={column name={$n_x$}, column type={r}}, %
    columns/nt/.style={column name={$n_t$}, column type={r}}, %
    columns/cgIterMean/.style={column name=\textbf{mean(CG)}, column type={r}, fixed, fixed zerofill, precision=2}, %
    columns/cgIterVar/.style={column name=\textbf{var(CG)}, column type={r}, fixed, fixed zerofill, precision=2}, %
    columns/l2Error/.style={column name={$\| u_{\varrho,h} -
        \overline{u}_s\|_{L^2(Q)}$}, column type={r}, sci, precision=3}, %
    columns/eoc/.style={column name=eoc, column type={r}, fixed, fixed zerofill, precision=2}, %
    ]{simulation/h22_3d/transformed_data.csv}
        \caption{Numerical results for a smooth target
      $\overline{u}_s \in H^{1,1/2}_{0;0,}(Q) \cap H^2(Q)$.}\label{tab:uS}

\end{table}  

\noindent  
As a second example we consider an anisotropic target
$\overline{u}_a \in H^{1,1/2}_{0;0,}(Q) \cap H^{2,1-\varepsilon}(Q)$,
$\varepsilon > 0$, i.e.,
\[
  \overline{u}_a(x, t) \, = \,
  \sqrt{t (T-t)} \, x_1 \, (1 - x_1) \, x_2 \, (1 - x_2) \,
  x_3 \, (1 - x_3) .
\]
In this case we have to use the parabolic scaling
$h_t = h_x^2$ in order to guarantee second order convergence, see
Theorem \ref{Thm error aniso}. The numerical results as given in
Table \ref{tab:uA:parabolicScaling} confirm this estimate.

\begin{table}[H]
    \centering
    \begin{filecontents*}{simulation/h21_3d_ps/transformed_data.csv}
dof,simulationTime,spatialProblemsTime,percentageSpatialSolve,nx,nt,eoc,l2Error,timePerDof
4,10,1,10.0,2,4,0.0,0.00221487,2500.0
432,60,28,46.666666666666664,4,16,0.402711,0.0016754,138.88888888888889
21952,902,658,72.9490022172949,8,64,0.995763,0.000840166,41.089650145772595
864000,25316,20507,81.00410807394533,16,256,1.57242,0.000282501,29.300925925925927
30505984,831581,686804,82.59015056861568,32,1024,1.85818,7.79202e-05,27.259602575022658
1024192512,26825636,22166773,82.63279573315614,64,4096,1.95486,2.00992e-05,26.191986062870182
\end{filecontents*}
\pgfplotstabletypeset[%
    col sep=comma, %
    columns={dof,nx,nt,l2Error,eoc, simulationTime, timePerDof}, %
    every head row/.style={before row=\toprule,after row=\midrule}, %
      every last row/.style={after row=\bottomrule}, %
      columns/simulationTime/.style={column name={time [ms]}, column type={r}, fixed}, %
      columns/timePerDof/.style={column name={time/dof [ns]}, column type={r}, fixed}, %
      columns/dof/.style={column name={DoF}, column type={r}, fixed}, %
      columns/nx/.style={column name={$n_x$}, column type={r}}, %
    columns/nt/.style={column name={$n_t$}, column type={r}, fixed}, %
    columns/cgIterMean/.style={column name=\textbf{mean(CG)}, column type={r}, fixed, fixed zerofill, precision=2}, %
    columns/cgIterVar/.style={column name=\textbf{var(CG)}, column type={r}, fixed, fixed zerofill, precision=2}, %
    columns/l2Error/.style={column name={
        $\| u_{\varrho,h} - \overline{u}_a\|_{L^2(Q)}$}, column type={r}, sci, precision=3}, %
    columns/eoc/.style={column name=eoc, column type={r}, fixed, fixed zerofill, precision=2}, %
    ]{simulation/h21_3d_ps/transformed_data.csv}
    \caption{Numerical results for an anisotropic target
      $\overline{u}_a \in H^{1,1/2}_{0;0,}(Q) \cap H^{2,1-\varepsilon}(Q)$,
      $\varepsilon > 0$, when using the parabolic scaling $h_t=h_x^2$.}
    \label{tab:uA:parabolicScaling}
  \end{table}  

\noindent
It is worth to mention that in this case, the use of $h_t=h_x$ is not
sufficient to obtain optimal convergence. Following the proof of Lemma
\ref{Lemma approximation aniso} we then conclude only linear convergence,
as confirmed by the numerical results given in Table \ref{tab:uA}.

\begin{table}[H]
    \centering
    \begin{filecontents*}{simulation/h21_3d/transformed_data.csv}
dof,simulationTime,spatialProblemsTime,percentageSpatialSolve,nx,nt,eoc,l2Error,timePerDof
2,4,0,0.0,2,2,0.0,0.0020859,2000.0
108,12,5,41.666666666666664,4,4,0.448338,0.00152873,111.11111111111111
2744,112,83,74.10714285714286,8,8,0.876656,0.000832589,40.816326530612244
54000,1628,1317,80.8968058968059,16,16,1.11806,0.000383585,30.14814814814815
953312,25553,20910,81.82992212264705,32,32,1.04226,0.000186256,26.804445973616193
16003008,421127,348018,82.63967876673783,64,64,0.982874,9.424e-05,26.315490187844684
262193024,6682653,5520303,82.60645884201978,128,128,0.979323,4.78002e-05,25.487531659118435
\end{filecontents*}
\pgfplotstabletypeset[%
    col sep=comma, %
    columns={dof,nx,nt,l2Error,eoc, simulationTime, timePerDof}, %
    every head row/.style={before row=\toprule,after row=\midrule}, %
      every last row/.style={after row=\bottomrule}, %
      columns/simulationTime/.style={column name={time [ms]}, column type={r}, fixed}, %
      columns/timePerDof/.style={column name={time/dof [ns]}, column type={r}, fixed}, %
      columns/dof/.style={column name={DoF}, column type={r}, fixed}, %
      columns/nx/.style={column name={$n_x$}, column type={r}}, %
    columns/nt/.style={column name={$n_t$}, column type={r}}, %
    columns/cgIterMean/.style={column name=\textbf{mean(CG)}, column type={r}, fixed, fixed zerofill, precision=2}, %
    columns/cgIterVar/.style={column name=\textbf{var(CG)}, column type={r}, fixed, fixed zerofill, precision=2}, %
    columns/l2Error/.style={column name={
        $\| u_{\varrho,h} - \overline{u}_a \|_{L^2(Q)}$}, column type={r}, sci,
    precision=3}, %
    columns/eoc/.style={column name=eoc, column type={r}, fixed, fixed zerofill, precision=2}, %
    ]{simulation/h21_3d/transformed_data.csv}

    \caption{Numerical results for an anisotropic target
      $\overline{u}_a \in H^{1,1/2}_{0;0,}(Q) \cap H^{2,1-\varepsilon}(Q)$,
      $\varepsilon > 0$ in the case of a uniform mesh with $h_t=h_x$.}
    \label{tab:uA}
\end{table}  

\noindent
The third example is a discontinuous target
$ \overline{u}_d \in H^{1/2-\varepsilon}(Q)$, $\varepsilon>0$,
\[
  \overline{u}_D(x, t) \, = \,
  \begin{cases}
    1 & \text{for } x \in (0.25, 0.75)^3, \\
    0 & \text{else}.
    \end{cases}
\]
In this case we have the error estimate \eqref{final error s} for all
$s < 1/2$. This convergence is again confirmed by the numerical results,
see Table \ref{tab:uD}.

\begin{table}[H]
    \centering
    \begin{filecontents*}{simulation/l2_3d/transformed_data.csv}
dof,simulationTime,spatialProblemsTime,percentageSpatialSolve,nx,nt,eoc,l2Error,timePerDof
2,6,0,0.0,2,2,0.0,0.458641,3000.0
108,18,8,44.44444444444444,4,4,1.10506,0.213215,166.66666666666666
2744,153,114,74.50980392156863,8,8,0.348337,0.167478,55.75801749271137
54000,1635,1329,81.28440366972477,16,16,0.459194,0.121822,30.27777777777778
953312,25629,21196,82.70318779507589,32,32,0.486596,0.0869453,26.88416803732671
16003008,419824,349444,83.23583215823774,64,64,0.49389,0.0617406,26.2340679952169
262193024,6756770,5636113,83.41430890795453,128,128,0.49679,0.0437544,25.770212711685264
\end{filecontents*}
\pgfplotstabletypeset[%
    col sep=comma, %
    columns={dof,nx,nt,l2Error,eoc, simulationTime, timePerDof}, %
    every head row/.style={before row=\toprule,after row=\midrule}, %
      every last row/.style={after row=\bottomrule}, %
      columns/simulationTime/.style={column name={time [ms]}, column type={r}, fixed}, %
      columns/timePerDof/.style={column name={time/dof [ns]}, column type={r}, fixed}, %
      columns/dof/.style={column name={DoF}, column type={r}, fixed}, %
      columns/nx/.style={column name={$n_x$}, column type={r}}, %
    columns/nt/.style={column name={$n_t$}, column type={r}}, %
    columns/cgIterMean/.style={column name=\textbf{mean(CG)}, column type={r}, fixed, fixed zerofill, precision=2}, %
    columns/cgIterVar/.style={column name=\textbf{var(CG)}, column type={r}, fixed, fixed zerofill, precision=2}, %
    columns/l2Error/.style={column name={
        $\| u_{\varrho,h} - \overline{u}_d\|_{L^2(Q)}$}, column type={r},
      sci,precision=3}, %
    columns/eoc/.style={column name=eoc, column type={r}, fixed, fixed zerofill, precision=2}, %
    ]{simulation/l2_3d/transformed_data.csv}
    \caption{Numerical results for the discontinuous target
      $\overline{u}_d \in H^{1/2-\varepsilon}(Q)$, $\varepsilon > 0$.}
    \label{tab:uD}
\end{table}  

\noindent
Finally, we consider the turning wave example in two space dimensions,
i.e., $\Omega = (0,1)^2$, subject to the heat equation with
a nonlinear (cubic) reaction term
$R(u) = u(u+1)(u-1/4)$, see \cite{Casas:2013,LSTY_SISC:2021}.
In order to have homogeneous initial and boundary conditions, we define
\begin{align*}
    \overline{u}_w(x, t) &=  \left(1 + \exp\left(\frac{\cos(g(t))\left(\frac{70}{3} - 70x_1\right) + \sin(g(t))\left(\frac{70}{3} - 70x_2\right)}{\sqrt{2}}\right)\right)^{-1} \\
    &+ \left(1 + \exp\left(\frac{\cos(g(t))\left(70x_1 - \frac{140}{3}\right) + \sin(g(t))\left(70x_2 - \frac{140}{3}\right)}{\sqrt{2}}\right)\right)^{-1} - 1
\end{align*}
for $x \in (1/8,7/8)^2$, and $t \in (1/8,1)$, and zero else, and where
$g(t) = \frac{2\pi}{3} \min \{ \frac{3}{4}, t \}$. Since the target is
discontinuous, $u_w \in H^{1/2-\varepsilon}(Q)$, $\varepsilon > 0$, follows,
and we can use the error estimate \eqref{final error s} for
$s < 1/2$, see Table \ref{tab:uT} for the related numerical results.

  \begin{table}[H]
    \centering
    \begin{filecontents*}{simulation/turning_wave_2d/transformed_data.csv}
dof,simulationTime,spatialProblemsTime,percentageSpatialSolve,nx,nt,eoc,l2Error,timePerDof
2,44,6,13.636363636363637,2,2,0.0,0.421184,22000.0
36,7,0,0.0,4,4,0.483014,0.301349,194.44444444444446
392,16,1,6.25,8,8,0.320045,0.241394,40.816326530612244
3600,49,11,22.448979591836736,16,16,0.662863,0.15247,13.61111111111111
30752,193,68,35.233160621761655,32,32,0.648301,0.0972808,6.2760145681581685
254016,1213,523,43.11624072547403,64,64,0.5948,0.0644132,4.775289745527841
2064512,9087,4141,45.57059535600308,128,128,0.532342,0.0445373,4.401524428048856
16646400,73386,34374,46.839996729621454,256,256,0.506822,0.0313441,4.408520761245675
133693952,633011,318031,50.24099107282496,512,512,0.501104,0.0221466,4.734776633725361
\end{filecontents*}
\pgfplotstabletypeset[%
    col sep=comma, %
    columns={dof,nx,nt,l2Error,eoc, simulationTime, timePerDof}, %
    every head row/.style={before row=\toprule,after row=\midrule}, %
      every last row/.style={after row=\bottomrule}, %
      columns/simulationTime/.style={column name={time [ms]}, column type={r}, fixed}, %
      columns/timePerDof/.style={column name={time/dof [ns]}, column type={r}, fixed}, %
      columns/dof/.style={column name={DoF}, column type={r}, fixed}, %
      columns/nx/.style={column name={$n_x$}, column type={r}}, %
    columns/nt/.style={column name={$n_t$}, column type={r}}, %
    columns/cgIterMean/.style={column name=\textbf{mean(CG)}, column type={r}, fixed, fixed zerofill, precision=2}, %
    columns/cgIterVar/.style={column name=\textbf{var(CG)}, column type={r}, fixed, fixed zerofill, precision=2}, %
    columns/l2Error/.style={column name=
      $\| u_{\varrho,h} - \overline{u}_w\|_{L^2(Q)}$, column type={r},
      sci, precision=3}, %
    columns/eoc/.style={column name=eoc, column type={r}, fixed, fixed zerofill, precision=2}, %
    ]{simulation/turning_wave_2d/transformed_data.csv}
    \caption{Numerical results for the turning wave $\overline{u}_w \in
      H^{1/2-\varepsilon}(Q)$, $\varepsilon > 0$.}
    \label{tab:uT}
\end{table}

\begin{figure}[H]
  \centering
  \begin{tabular}{ccc}
    \includegraphics[width=0.3\linewidth]{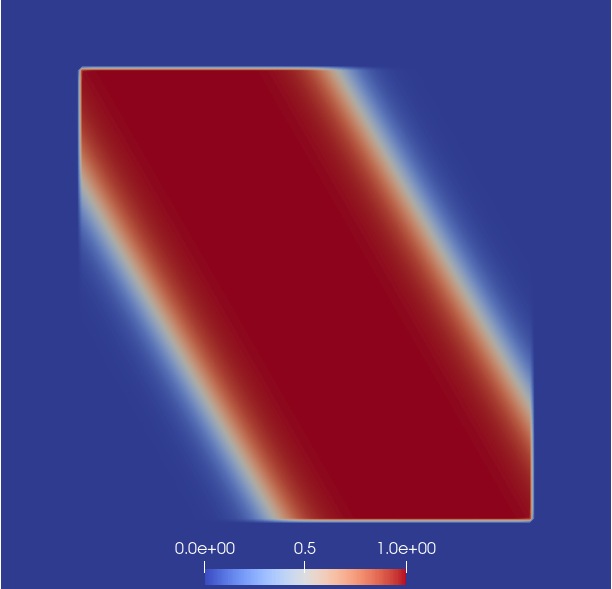} & \includegraphics[width=0.3\linewidth]{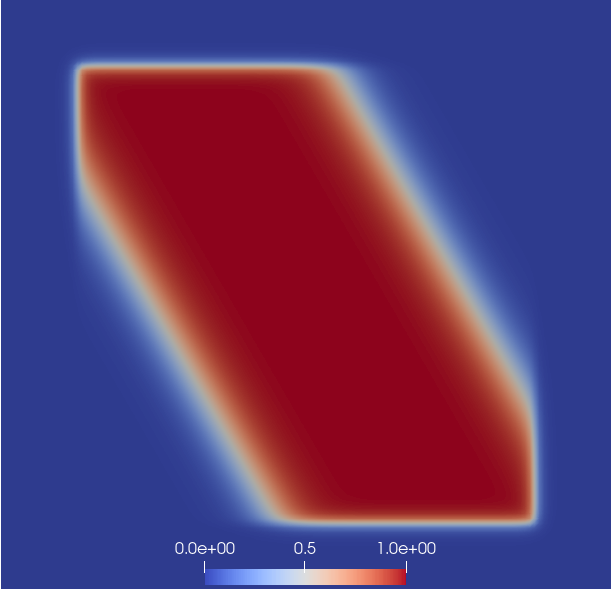} & \includegraphics[width=0.3\linewidth]{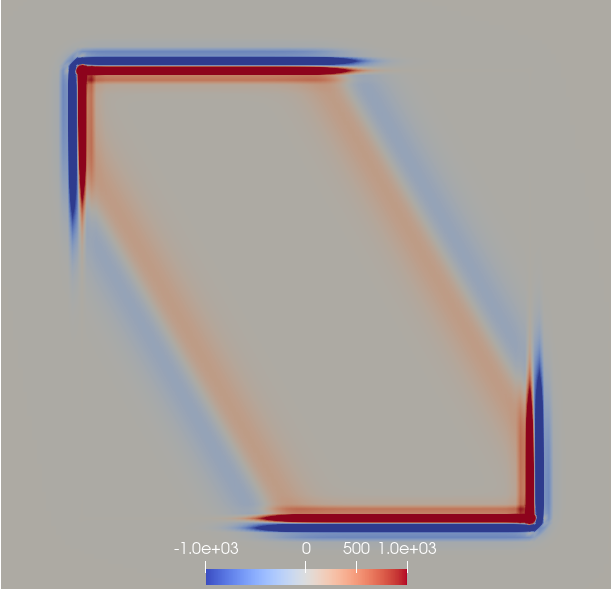} \\
    \includegraphics[width=0.3\linewidth]{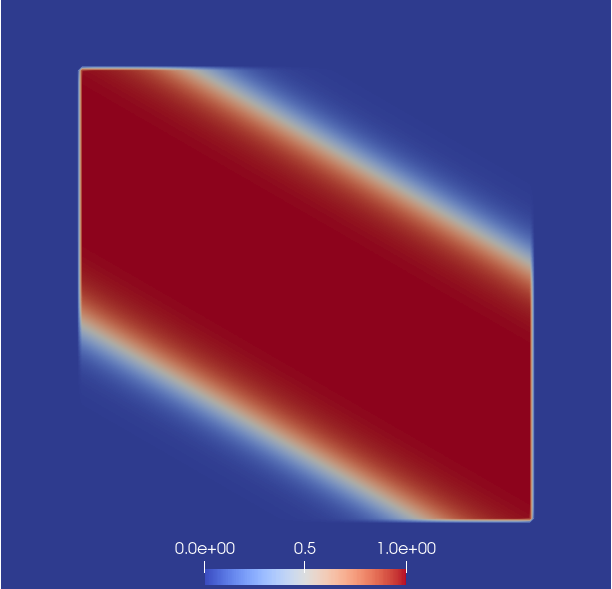} & \includegraphics[width=0.3\linewidth]{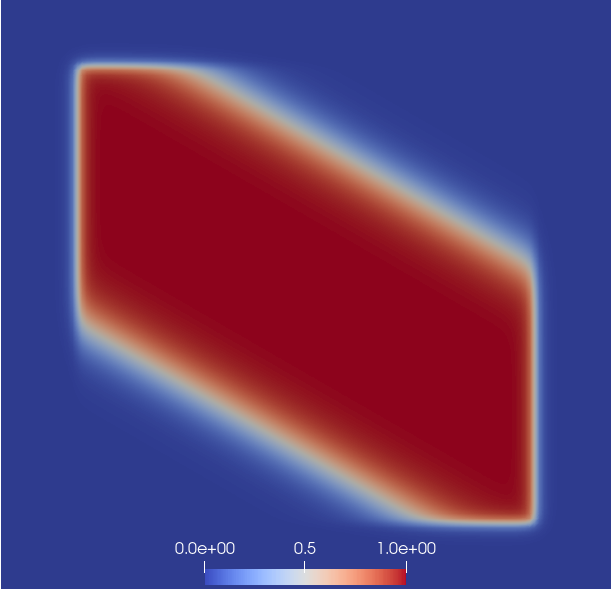} & \includegraphics[width=0.3\linewidth]{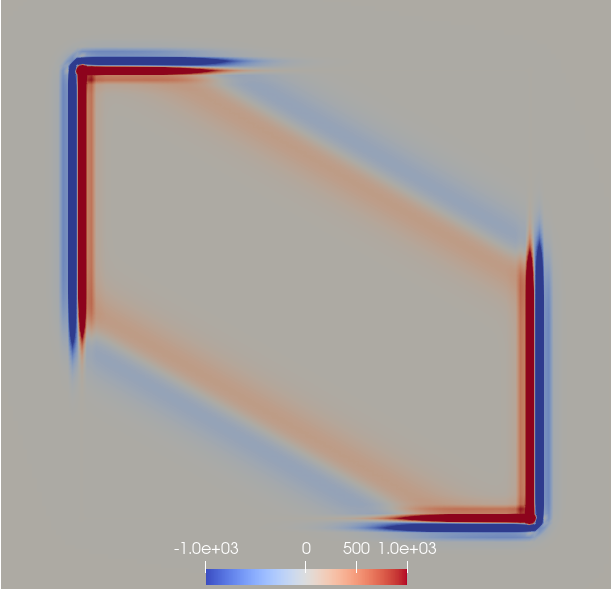} \\
    \includegraphics[width=0.3\linewidth]{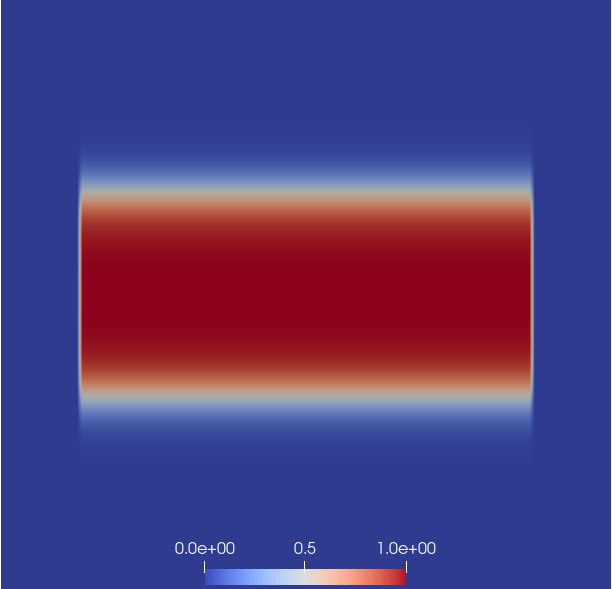} & \includegraphics[width=0.3\linewidth]{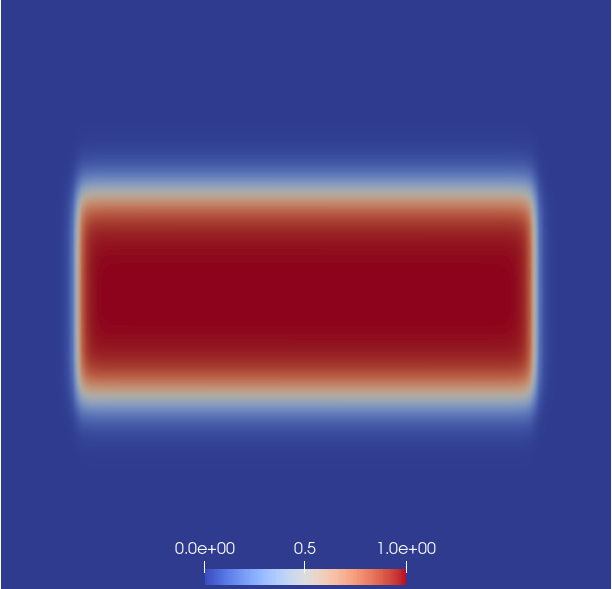} & \includegraphics[width=0.3\linewidth]{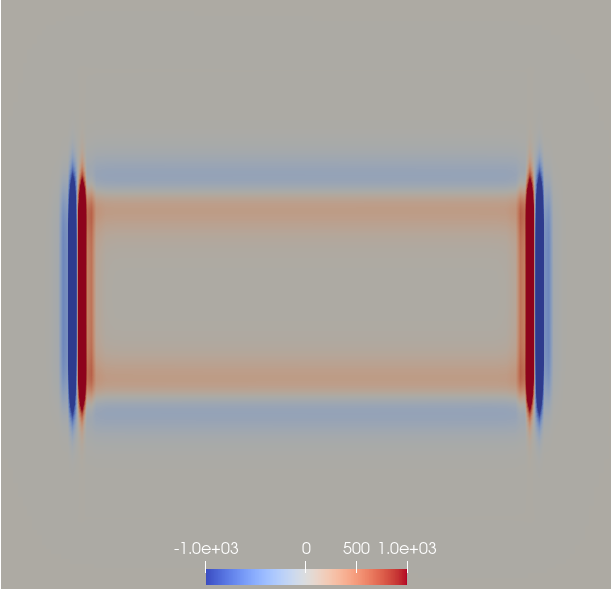} \\
  \end{tabular}

  \caption{Simulation results for the nonlinear turning wave example with
    $n_t=n_x=128$. 
    Target (left), state (middle), and control (right) at $t=0.25$ (top),
    $t=0.5$ (middle), and $t=0.75$ (bottom).}
\label{fig:turningWave}

\end{figure}

\noindent
Once we have computed an approximation $u_{\varrho,h}$ of the state
$u_\varrho \in H^{1,1/2}_{0;0,}(Q)$ we can recover the corresponding control
$z_\varrho = B u_\varrho \in [H^{1,1/2}_{0;,0}(Q)]^*$ via post processing,
see \cite{LLSY:adaptive} in the elliptic case. When using a least squares
approach, we have to solve a saddle point formulation to find
$(w,z_\varrho) \in H^{1,1/2}_{0;,0}(Q) \times [H^{1,1/2}_{0;,0}(Q)]^*$ such
that
\begin{equation}\label{mixed VF control}
  \langle A w , v \rangle_Q + \langle z , v \rangle_Q =
  \langle B u_\varrho , v \rangle_Q, \quad
  \langle w , \psi \rangle_Q = 0
\end{equation}
is satisfied for all $(v,\psi) \in
H^{1,1/2}_{0;,0}(Q) \times [H^{1,1/2}_{0;,0}(Q)]^*$, where
$A : H^{1,1/2}_{0;,0}(Q) \to [H^{1,1/2}_{0;,0}(Q)]^*$ is the Riesz
operator as defined in \eqref{Def A aniso}.
Using a stable discretization of \eqref{mixed VF control}, e.g., using
a piecewise linear continuous approximation $w_h$, and a piecewise
constant approximation $z_{\varrho, h}$ on a coarser mesh, and replacing
the state $u_\varrho$ by its approximation $u_{\varrho,h}$, we obtain
an approximate control. But here we apply a simpler approach.
From the optimality system we easily conclude $z_\varrho \in L^2(Q)$.
Hence we can compute $z_{\varrho,h} = Q_h B u_{\varrho,h}$ as $L^2$
projection on a suitable ansatz space. When using a piecewise (multi-)linear
approximation $u_{\varrho,h}$ as in this paper, we have to consider a piecewise
linear approximation for the control as well in order to apply integration
by parts to evaluate the right hand side. Alternatively, we can use
second order B splines in all spatial coordinates, and still linear ones
in time, to allow for a direct evaluation of $B u_{\varrho,h}$, and hence,
we can use piecewise constants to approximate the control.
Figure \ref{fig:turningWave} depicts the turning wave target
$\overline{u}_w$ as well as the corresponding state $u_{\varrho,h}$, and
the reconstructed piecewise linear control $z_{\varrho,h}$.

\noindent
In Figure \ref{fig:timings} we finally summarize the computing times for
all examples, confirming, up to logarithmic factors, optimal complexity
in all cases.

  \begin{figure}[H]
    \centering
    \begin{tikzpicture}
      \begin{loglogaxis}[
          xlabel={dof},
          ylabel={Cpu time [ms]},
          width=0.8\textwidth, % Adjust the width as needed
          height=0.6\textwidth, % Adjust the height as needed
          grid=both,
          minor tick num=5,
          legend pos=north west,
          legend cell align={left}
      ]
      \addplot table [x=dof, y=simulationTime, col sep=comma] {simulation/h22_3d/transformed_data.csv};
      \addlegendentry{$\overline{u}_s$}
      \addplot table [x=dof, y=simulationTime, col sep=comma] {simulation/h21_3d/transformed_data.csv};
      \addlegendentry{$\overline{u}_a$}
      \addplot table [x=dof, y=simulationTime, col sep=comma] {simulation/h21_3d_ps/transformed_data.csv};
      \addlegendentry{$\overline{u}_a$ (parabolic scaling)}
      \addplot table [x=dof, y=simulationTime, col sep=comma] {simulation/l2_3d/transformed_data.csv};
      \addlegendentry{$\overline{u}_d$}
      \addplot table [x=dof, y=simulationTime, col sep=comma] {simulation/turning_wave_2d/transformed_data.csv};
      \addlegendentry{$\overline{u}_w$}
      \end{loglogaxis}
  \end{tikzpicture}
    \caption{Simulation times for all cases.}
    \label{fig:timings}
\end{figure}

\begin{figure}[htbp!]
  \centering
  \begin{tabular}{c c }
    \includegraphics[width=0.45\textwidth]{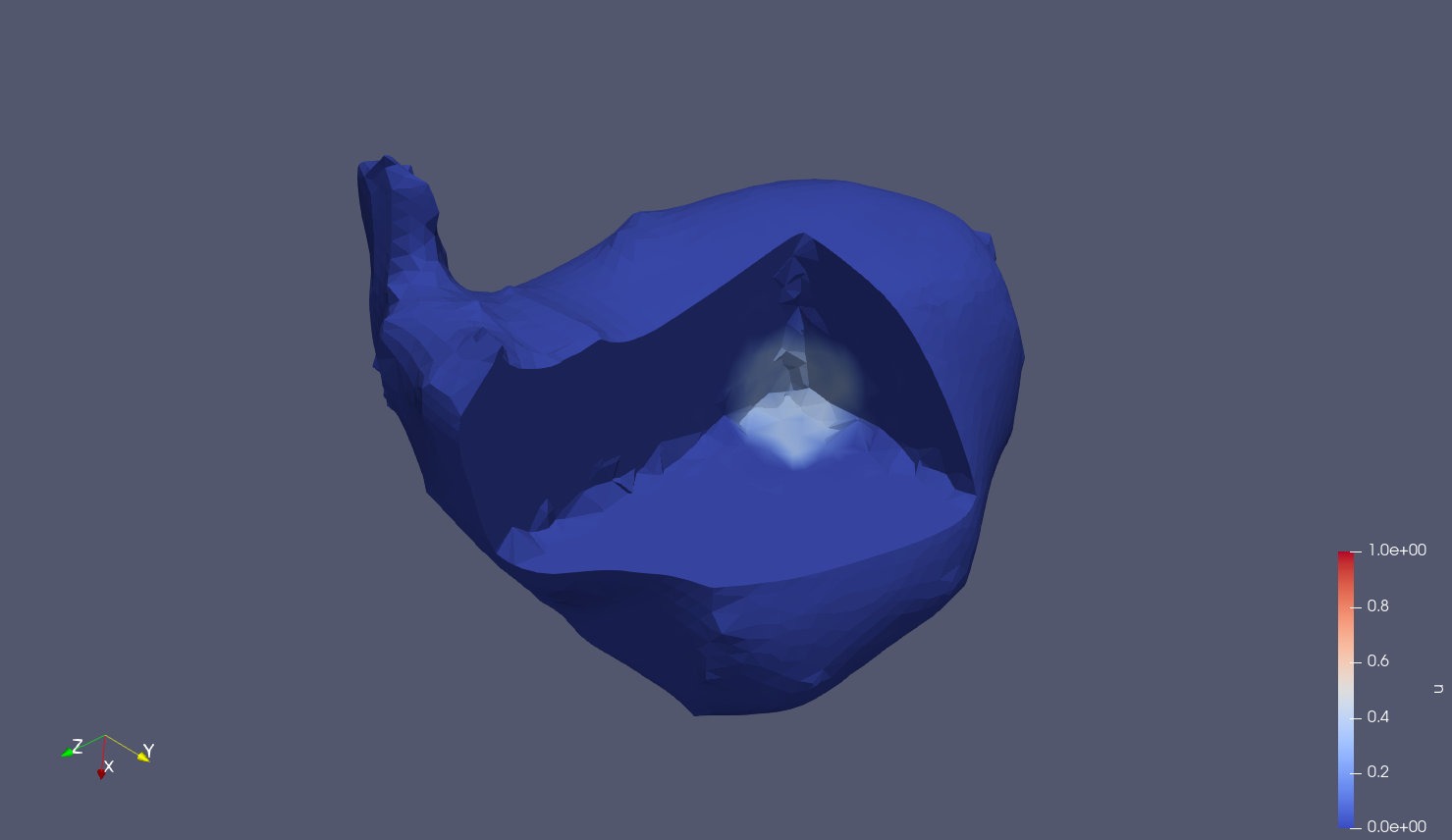} 
    & \includegraphics[width=0.45\textwidth]{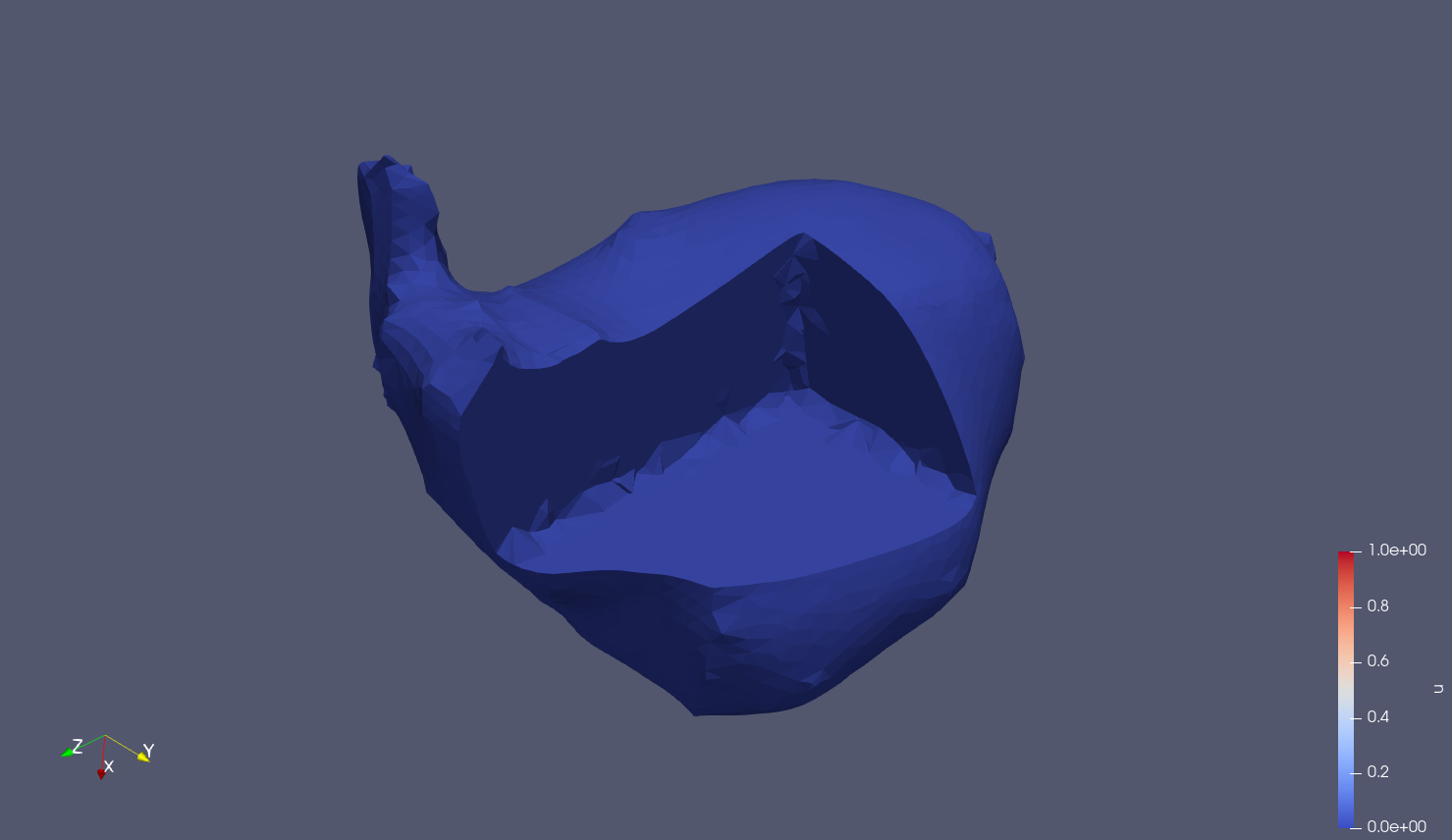} \\
    \includegraphics[width=0.45\textwidth]{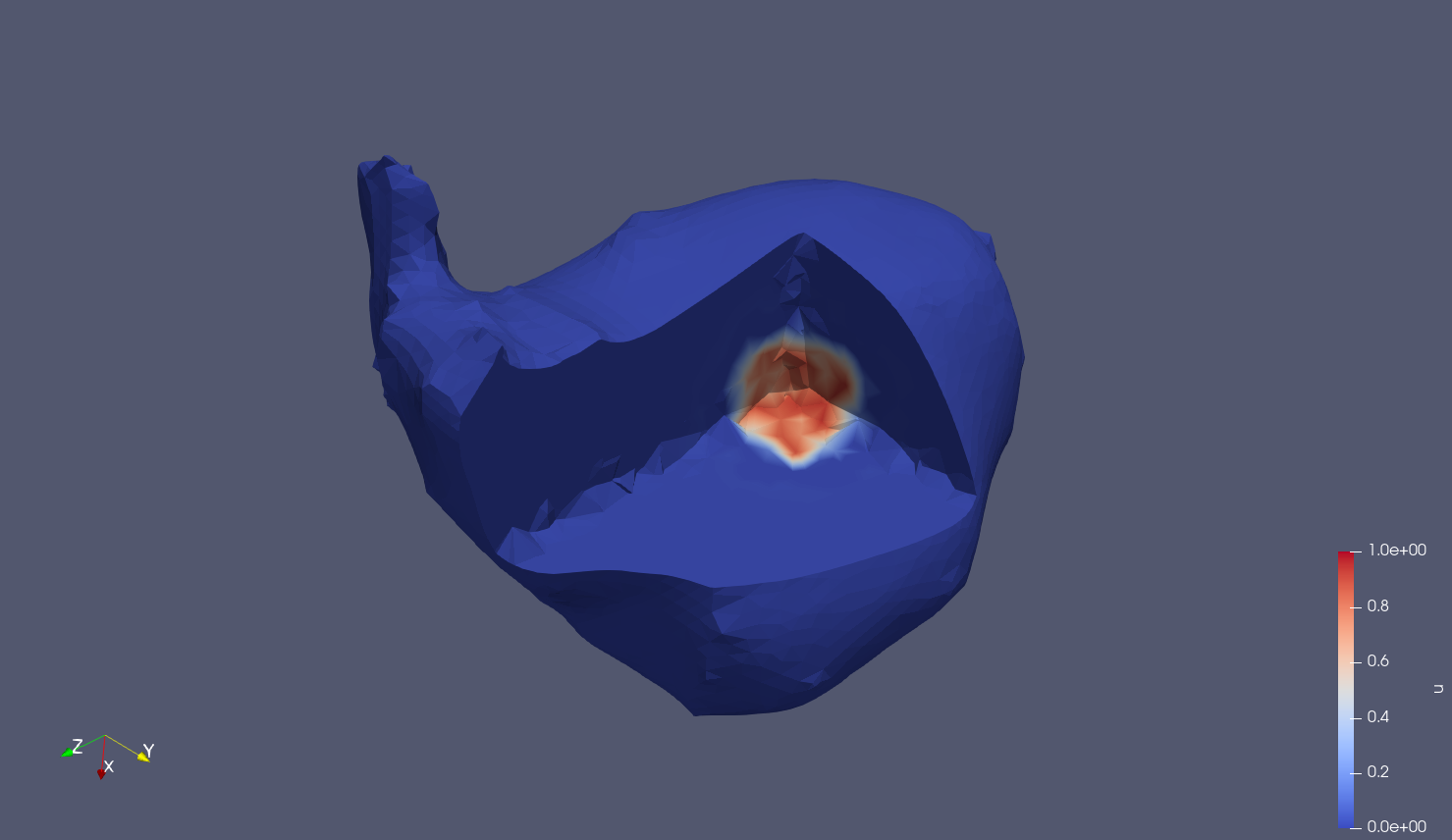} 
    & \includegraphics[width=0.45\textwidth]{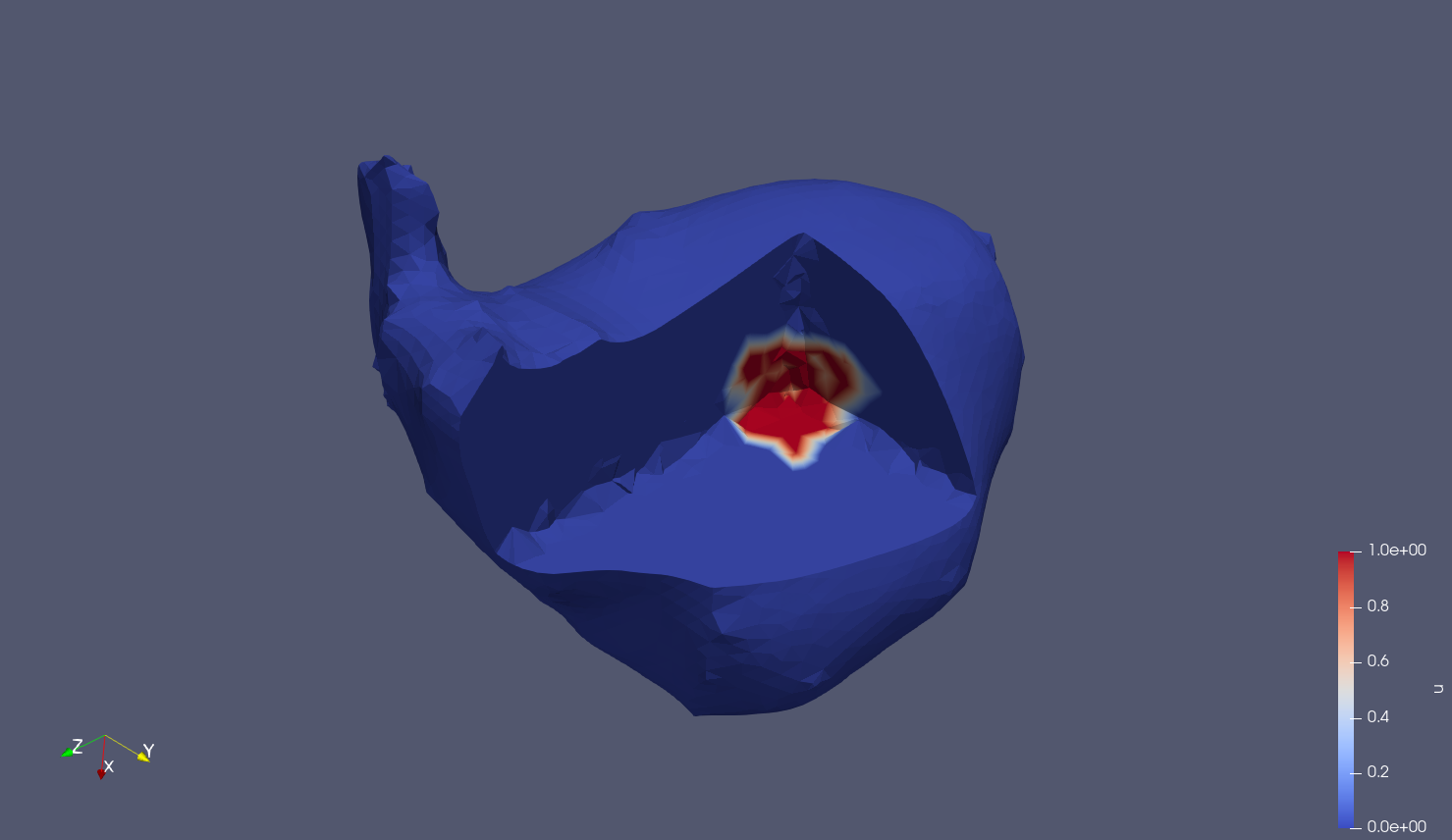} \\
    \includegraphics[width=0.45\textwidth]{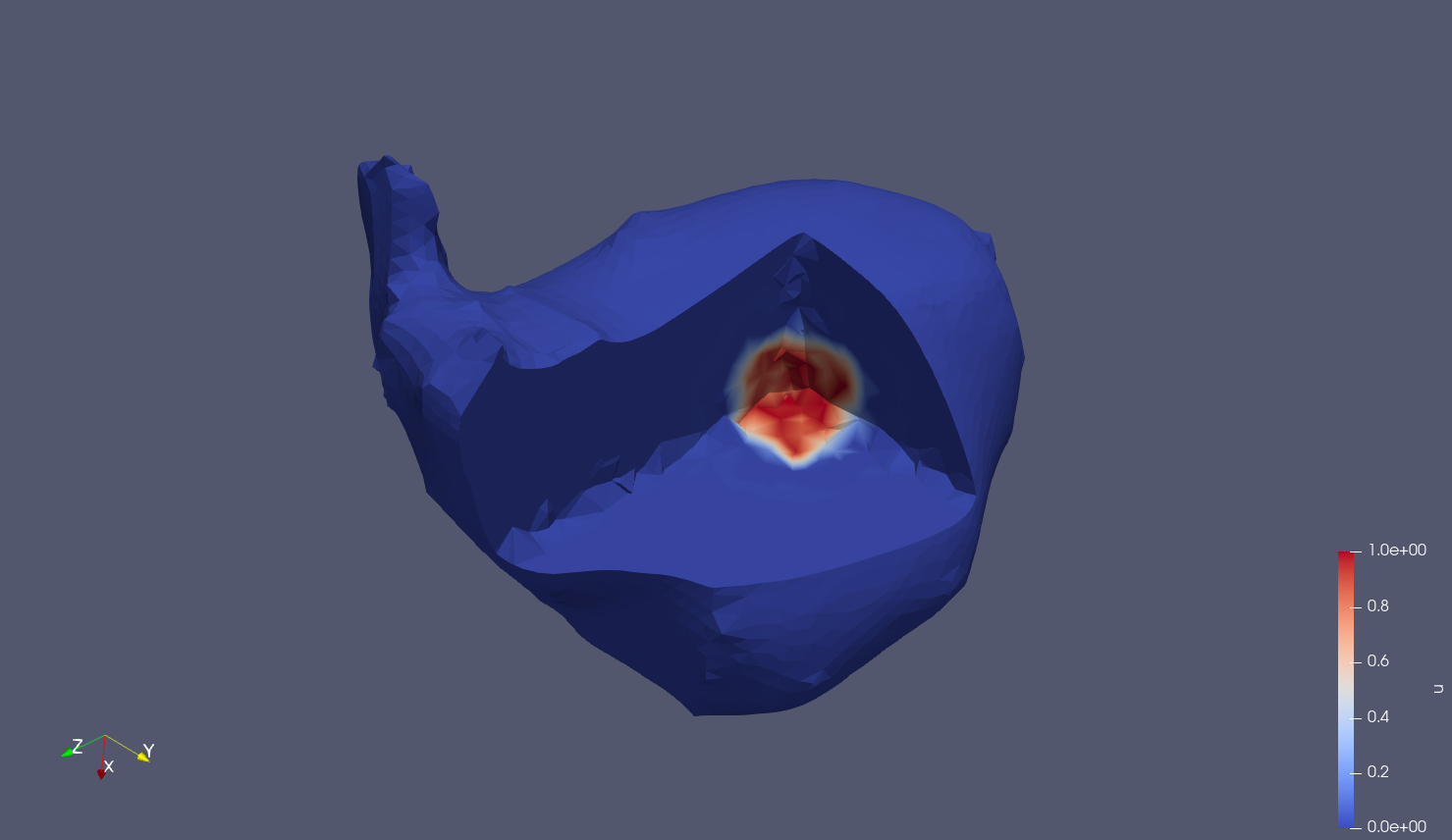} 
    & \includegraphics[width=0.45\textwidth]{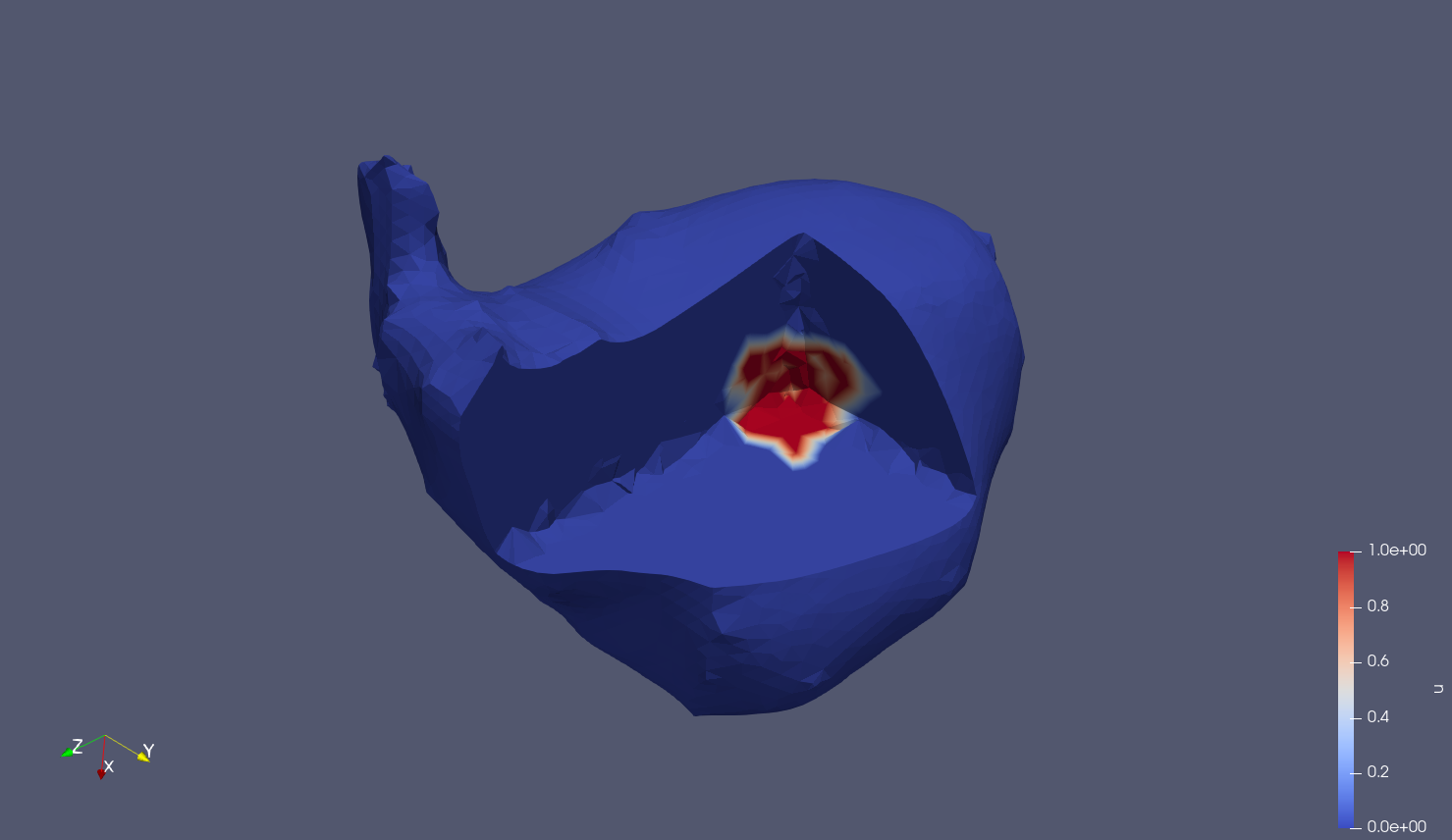}

  \end{tabular} 
  \caption{Gallbladder with discontinuous target function, \\
  left column: state $u_{\varrho,h}$, right column: target $\overline{u}$ for times $t=0.24$ (top), $t=0.3$ (middle), and $t=0.5$ (bottom).}
  \label{fig:gallbladder} 
\end{figure}

\noindent
As a last, more practical, example, we consider the heating of a gallbladder,
see Figure~\ref{fig:gallbladder}. We aim to approximate a constant
cube-shaped target function, that is turned on at $t=0.25$, i.e.
\begin{align*}
  \overline{u}(x,t) = \begin{cases}
    1 & \text{for } t>0.25 \text{ and } x \in C, \\
    0 & \text{else}.
    \end{cases}
\end{align*}
where $C$ is the spatial cubic region. Further we impose the condition, that the state must be smaller or equal to 1 at all times.
Moreover the state $u_\varrho$ is satisfying the heat equation with homogenuous boundary and initial conditions. Note, that the given target satisfies boundary conditions and limits, but is discontinuous both in space and time.  
The spatial mesh is generated from a surface mesh \cite{githubgallbladder}
using gmsh \cite{geuzaine2009gmsh}, contains $5582$ nodes and is completely
unstructured. The spatial finite element matrices and vectors are computed
using MFEM \cite{mfem}. The temporal mesh consists of $50$ equidistant temporal
elements on the interval $(0,1)$. The starting point for the semismooth Newton
method is the interpolation of the target for the state and the Lagrange
parameters are initialized with zero. As the focus of this paper is not on the
nonlinear solver, we used a rather simple stopping criterion, where we stop as
soon as subsequent iterates change less than $10^{-3}$ in the maximum norm and used aggressive damping of $0.01$. For details on the nonlinear solver for state constraints see \cite{löscher2024optimalcomplexitysolutionspacetime}. In the temporal development of state and target, depicted in Figure \ref{fig:gallbladder}, we can see that the state precedes the target a little and needs more time to reach the final state, which is reasonable due to the fact that the target is discontinuous.Moreover it never exceeds the value of 1.

